\documentclass[11pt]{article}
\usepackage{amsthm, amsmath, amssymb, amsfonts, url, booktabs, tikz, setspace, fancyhdr, amsbsy}
\usepackage[margin = 0.8in]{geometry}
\usepackage{hyperref, enumerate}
\usepackage{caption,float}
\usepackage{subcaption}
\usepackage{esint}
\usepackage{verbatim}
\usepackage{multirow}
\usepackage{kotex}
\usepackage{graphicx, wrapfig}

\newtheorem{theorem}{Theorem}[section]

\newtheorem{lemma}[theorem]{Lemma}

\theoremstyle{definition}

\newcommand{\R}{\mathbb{R}}

\newcommand{\dx}{\,\mathrm{d}x}
\newcommand{\dy}{\,\mathrm{d}y}
\newcommand{\dz}{\,\mathrm{d}z}
\newcommand{\xx}{\boldsymbol{x}}

\newcommand{\tth}{\boldsymbol{\theta}}

\numberwithin{equation}{section}

\begin{document}

\title{VS-PINN: A fast and efficient training of physics-informed neural networks using variable-scaling methods for solving PDEs with stiff behavior}

\author{Seungchan Ko\thanks{Department of Mathematics, Inha University, Incheon, Republic of Korea. Email: \tt{scko@inha.ac.kr}}
~and ~Sang Hyeon Park\thanks{Department of Mathematics, Inha University, Incheon, Republic of Korea. Email: \tt{12181868@inha.edu}}}

\date{~}

\maketitle

~\vspace{-1.5cm}


\begin{abstract}
    Physics-informed neural networks (PINNs) have recently emerged as a promising way to compute the solutions of partial differential equations (PDEs) using deep neural networks. However, despite their significant success in various fields, it remains unclear in many aspects how to effectively train PINNs if the solutions of PDEs exhibit stiff behaviors or high frequencies. In this paper, we propose a new method for training PINNs using variable-scaling techniques. This method is simple and it can be applied to a wide range of problems including PDEs with rapidly-varying solutions. Throughout various numerical experiments, we will demonstrate the effectiveness of the proposed method for these problems and confirm that it can significantly improve the training efficiency and performance of PINNs. Furthermore, based on the analysis of the neural tangent kernel (NTK), we will provide theoretical evidence for this phenomenon and show that our methods can indeed improve the performance of PINNs.
\end{abstract}


\smallskip



\section{Introduction}
Bridging the gap between classical scientific computing and machine learning, the emerging field called scientific machine learning has introduced a completely different framework to compute the solutions of partial differential equations (PDEs). The forefront of this evolution lies within the area of physics-informed neural networks (PINNs) \cite{pinn01, pinn02}. Based on the universal approximation property of deep neural networks, the PINNs approximate the solutions of PDEs by incorporating underlying physics into deep neural networks. Due to their simplicity and flexibility in handling a wide range of physical problems involving PDEs, PINNs have recently gained great attention and have been applied to various fields in computational science: bio-medical science \cite{pinn_bio_1, pinn_bio_2}, fluids mechanics \cite{pinn_fluid_1, pinn_fluid_2, pinn_fluid_3, pinn_fluid_4, NS_comp}, uncertainty quantification \cite{pinn_uq_1, pinn_uq_2, pinn_uq_3} and meta-material design \cite{pinn_meta_1, pinn_meta_2}. Moreover, since the PINNs utilize randomly-selected collocation points as training samples in the spatio-temporal domain, the PINNs are available for high-dimensional PDEs \cite{pinn_hd_1, pinn_hd_2}, on the domains with complex geometries \cite{PINN_field_1, PINN_field_2, PINN_field_3, PINN_domain}.

However, despite the significant empirical success of PINNs, we only have limited knowledge about the behavior of these constrained neural networks during their training process, and the training of PINNs often fails. In particular, since neural networks typically assume a smooth prior, it is often challenging to train PINNs to learn a solution with a sharp transition which poses significant obstacles for the model prediction. For example, some recent studies have illustrated that training PINNs with fully-connected neural networks usually suffers from so-called {\textit{spectral bias}} or {\textit{F-principle}} meaning that it is difficult for PINNs to learn functions with high frequencies \cite{bias_1, bias_2, bias_3}. Furthermore, it has recently been observed that the PINNs have some problems regarding convergence and accuracy when the target functions have sharp spatial transitions or fast temporal evolution \cite{sharp_ex}. Some previous works investigated the ways to handle these issues. For example, \cite{pre_work_1} theoretically analyzed the training process of PINNs based on the recently developed neural tangent kernel (NTK) \cite{NTK}, which explores the connection between deep neural networks and the kernel regression, elucidating the training dynamics. From the intuition obtained in the analysis, the authors proposed an NTK-guided gradient-descent algorithm to efficiently train PINNs without spectral-bias pathology and achieved remarkable improvements on various benchmark problems. In order to handle stiff PDEs, the authors of \cite{pre_work_2} developed Self-Adaptive PINNs, a novel paradigm for training PINNs using trainable weights based on the attention mechanism. On the other hand, for the singular perturbation problem which is another area in stiff PDEs, \cite{pre_work_3} developed a semi-analytic PINN
method, based on the so-called corrector functions obtained from the boundary layer analysis. 

In contrast to these methods, the present paper aims to provide a completely different method for tackling the aforementioned problems. In particular, we will propose a novel training method for PINNs using variable-scaling, which we have named VS-PINNs (variable-scaling physics-informed neural networks). Scaling is a common technique in the analysis of PDEs, but to the best of our knowledge, the idea has never been applied to train PINNs. Throughout various numerical experiments, we will show that this method can greatly improve the performance of PINNs and it is very effective in solving stiff problems. As it will be made clear soon, our method is simple and no significant extra computational costs are incurred. 
 Furthermore, it also offers computational flexibility, resulting in significant advancements in accuracy and learning efficiency for various computational problems. 

The rest of the paper is organized as follows. In Section \ref{sec:prelim}, we will collect some preliminary materials including brief descriptions of the PINNs and the NTK, which will be considered throughout the paper. In Section \ref{sec:method}, we shall describe the method that we mainly propose in this paper. And then in Section \ref{sec:exp}, we will perform comprehensive experiments to evaluate the performance of the proposed method. Throughout the numerical simulation of a wide range of PDEs with stiff behavior, high frequency and even for nonlinear problems, we will illustrate the effectiveness of the method we propose. Once we empirically confirm the better performance of the new method, in Section \ref{sec:thoery}, we will also rigorously investigate the proposed method by using the NTK. By examining the eigenvalues of the NTK of VS-PINNs, we will provide some theoretical evidence on why the proposed method can indeed improve the performance of PINNs. Finally in Section \ref{sec:conclude}, we will make a concluding remark, addressing the limitation and further improvement of the method, as well as the future research direction regarding the VS-PINNs.


\section{Preliminary}\label{sec:prelim}

\subsection{Physics-Informed Neural Networks}

\begin{figure}[h]
    \centering    \includegraphics[width=0.95\textwidth]{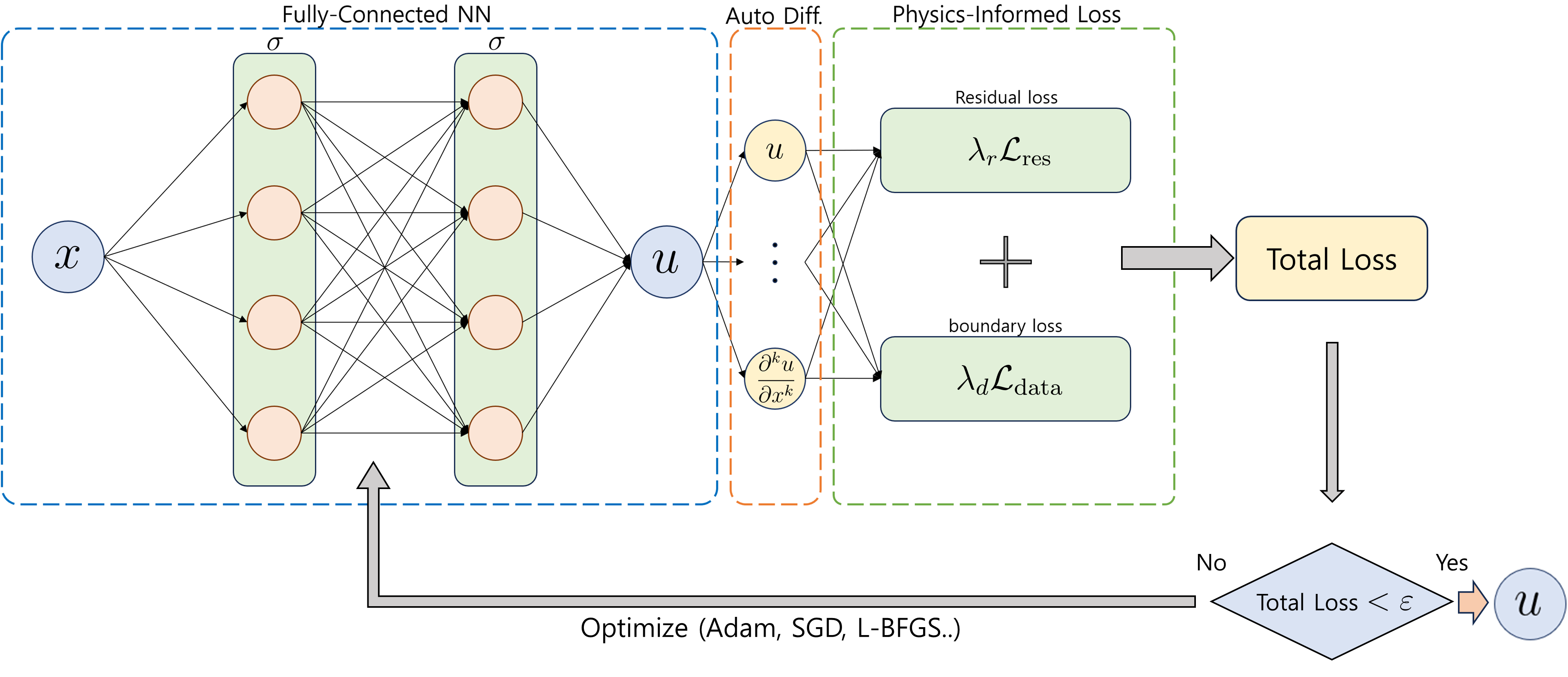}
    \caption{Schematic illustration of the physics-informed neural network. The left part visualizes a standard feed-forward neural
network parameterized by $\theta$, while the right part imposes the given physical laws to the neural network.}
    \label{F1}
\end{figure}


In this section, we will briefly review the concept of physics-informed neural networks (PINN) \cite{pinn01}. For the description of the method, let us consider the following problem on a bounded open domain $\Omega\subset\R^d$:
\begin{equation}\label{main_eq}
\begin{aligned}
\mathcal{D}[u](\xx)&=f(\xx), &&\boldsymbol{x}\in\Omega,\\
u(\boldsymbol{x})&=g(\boldsymbol{x}),&&\boldsymbol{x}\in\partial\Omega,
\end{aligned}
\end{equation}
where $\mathcal{D}$ is a differential operator and $u:\overline{\Omega}\rightarrow\R$ is the unknown function with the variable $\xx=(x_1,x_2,\cdots,x_d)\in\Omega$. For the sake of simplicity, for time-dependent problems, we regard the temporal variable $t$ as an additional coordinate of $\xx$. Hence $\Omega$ may denote the spatio-temporal domain, and the initial condition can be regarded as a special type of boundary data on $\partial\Omega$.
The main idea of the PINN is to approximate the solution $u(\xx,t)$ with a deep neural network $u(\xx,t;\tth)$, where $\tth$ is a collection of all parameters of the neural network. To train the model, we minimize the so-called {\textit{physics-informed loss function}} consisting of the residual of the given PDE in the domain and the differences of the neural network with the imposed conditions on the boundary. More precisely,
the neural network is trained by minimizing the total loss 
\begin{equation}\label{total_loss}
    \mathcal{L}_{\rm total}=\lambda_r\mathcal{L}_{\rm res}+\lambda_d\mathcal{L}_{\rm data},
\end{equation}
where each loss function is defined as follows:
\begin{align}
\mathcal{L}_{\rm res}&=\frac{1}{N_r}\sum_{i=1}^{N_r}\left|\mathcal{D}[u](\boldsymbol{x}_r^i;\tth)-f(\boldsymbol{x}_r^i)\right|^2,\label{res_loss}\\
\mathcal{L}_{\rm data}&=\frac{1}{N_b}\sum_{j=1}^{N_b}\left|u(\boldsymbol{x}_b^j;\tth)-g(\boldsymbol{x}_b^j)\right|^2\label{bdry_loss}.
\end{align}
Here $N_r$, $N_b$ denote the batch sizes, $\lambda_r$, $\lambda_b$ are the weights for the residual and boundary data respectively. Typically, $\boldsymbol{x}_r^i$ and $\boldsymbol{x}_b^j$ are randomly sampled from the uniform distributions $\mathcal{U}(\Omega)$ and $\mathcal{U}(\partial\Omega)$ respectively, and all the derivatives involved in the loss function can be computed accurately and quickly using reverse-mode {\textit{automatic differentiation}} \cite{auto_diff_1, auto_diff_2}.
This ensures that information on initial/boundary conditions propagates the interior domain while satisfying the physics laws represented as PDEs, and is the main reason why we call it a ``physics-informed'' neural network. If we have any sensor data for training, in other words, the ground-truth data $(\xx^i_s,y^i_s)$ with $y^i_s=u(\xx_s^i)$ for all $i=1,2,\cdots,N_s$, we can further add the so-called supervised loss $\mathcal{L}_{\rm sup}=\sum^{N_s}_{i=1}|u(\xx^i_s;\tth)-y^i_s|^2$ to the total loss, which is known to improve the performance of the model. However, in this paper, we will not consider this type of loss function, and propose a novel method which can train the PINN efficiently without the use of the training data. The schematic illustration of the PINN is depicted in Figure \ref{F1}.


\subsection{Neural Tangent Kernel}
In the seminar work \cite{NTK}, the authors proposed a novel approach, the so-called neural tangent kernel (NTK), that opened up new avenues in the analysis of training neural networks. The authors showed that under a certain parameter initialization, as the width of the layer increases infinitely, the kernel becomes constant during training so that training neural networks is equivalent to performing a deterministic kernel regression. In other words, under appropriate conditions, analysis of training using gradient descent can be regarded as a kernel regression analysis. In this section, we will briefly introduce the theory of NTK, which will be used in the later section to theoretically investigate the method proposed in this paper.

Let us consider the scenario where we need to minimize the square loss function on a dataset $(\xx_i,y_i)^m_{i=1}$: $\frac{1}{2}\sum^m_{i=1}(y_i-u(\xx_i;\boldsymbol{\theta}))^2$, where $u$ is parametrized by $\boldsymbol{\theta}$. A continuous version of gradient descent (gradient flow) can be written as 
\begin{equation}\label{2_2_1}
    \frac{{\rm d}\boldsymbol{\theta}(t)}{{\rm d}t}=-\nabla_{\tth}\mathcal{L}(\boldsymbol{\theta}(t))=-\sum_{i=1}^m(y_i-u(\boldsymbol{x}_i;\boldsymbol{\theta}(t)))\nabla_{\tth}u(\boldsymbol{x}_i;\boldsymbol{\theta}(t)).
\end{equation}
Note further that
\begin{equation}\label{2_2_2}
    \frac{{\rm d}u(\boldsymbol{x};\boldsymbol{\theta}(t))}{{\rm d}t}=
    \bigg[{\frac{{\rm d}u(\boldsymbol{x};\boldsymbol{\theta}(t))}{{\rm d}\boldsymbol{\theta}}}\bigg]\cdot \bigg[\frac{{\rm d}\boldsymbol{\theta}(t)}{{\rm d}t}\bigg]=\sum_{i=1}^m(y_i-u(\boldsymbol{x}_i;\boldsymbol{\theta}(t)))[\nabla_ {\boldsymbol{\theta}}u(\boldsymbol{x}_i;\boldsymbol{\theta}(t))]^T[\nabla_{\tth}u(\boldsymbol{x};\boldsymbol{\theta}(t))].
\end{equation}
This motivates us to define a component of NTK in the following way:
\[\widehat{\boldsymbol{\Theta}}(\boldsymbol{x},\boldsymbol{x}')(t)=\left\langle\frac{\partial u(\boldsymbol{x};\boldsymbol{\theta}(t))}{\partial\tth},\frac{\partial u(\boldsymbol{x}';\boldsymbol{\theta}(t))}{\partial\tth}\right\rangle.
\]
Then over the dataset $(\xx_i,y_i)^m_{i=1}$, \eqref{2_2_2} can be rewritten as
\begin{equation}\label{2_2_3}
\frac{\rm d}{{\rm d} t}{\left[\begin{matrix}
    u(\boldsymbol{x}_1;\boldsymbol{\theta}(t)) \\
    u(\boldsymbol{x}_2;\boldsymbol{\theta}(t)) \\
    \vdots \\ 
    u(\boldsymbol{x}_m;\boldsymbol{\theta}(t))
\end{matrix}\right]}=-\widehat{\boldsymbol{\Theta}}(t)\left[\begin{matrix}
    y_1-u(\boldsymbol{x}_1;\boldsymbol{\theta}(t)) \\
    y_2-u(\boldsymbol{x}_2;\boldsymbol{\theta}(t)) \\
    \vdots \\
    y_3-u(\boldsymbol{x}_m;\boldsymbol{\theta}(t))
\end{matrix}\right],\end{equation}
where $\widehat{\boldsymbol{\Theta}}(t)\in\mathbb{R}^{m\times m}$ and the $(i, j)$-th component of $\widehat{\boldsymbol{\Theta}}(t)$ is 
\[\widehat{\boldsymbol{\Theta}}_{ij}(t) =\widehat{\boldsymbol{\Theta}}(\boldsymbol{x}_i,\boldsymbol{x}_j)=[\nabla_{\tth}u(\boldsymbol{x}_i;\boldsymbol{\theta}(t))]^T[\nabla_{\tth}u(\boldsymbol{x}_j;\boldsymbol{\theta}(t))]. \]
If we write $\boldsymbol{U}(t)=\{U_i(t)\}^m_{i=1}$ with  $U_i(t)=y_i-u(\boldsymbol{x}_i;\boldsymbol{\theta}(t))$, and suppose the NTK remains constant and positive definite during training (i.e. $\widehat{\boldsymbol{\Theta}}(t)=\widehat{\boldsymbol{\Theta}}(0)$ and all the eigenvalues of $\widehat{\boldsymbol{\Theta}}(0)$ are positive). Then from \eqref{2_2_3}, we have
\[
    \frac{\rm d}{{\rm d}t}\left(\frac{1}{2}\|\boldsymbol{U}(t)\|^2_2\right)=-\boldsymbol{U}^T(t)\widehat{\boldsymbol{\Theta}}(0)\boldsymbol{U}(t)\leq-\frac{\lambda_0}{2}\|\boldsymbol{U}(t)\|^2_2,
\]
where $\lambda_0$ is the smallest eigenvalue of $\widehat{\boldsymbol{\Theta}}(0)$, which leads us to the convergence result
\begin{equation}\label{2_2_4}
    \|\boldsymbol{U}(t)\|^2_2\leq e^{-\lambda_0t}\|\boldsymbol{U}(0)\|^2_2.
\end{equation}

Next, we shall discuss the assumption that the NTK remains constant during the training process, and we will briefly review the relevant theory. To do this, let us consider an $L$-layer fully-connected neural network $f(\xx)=g_L(\xx)$ defined recursively as
\[
    g_{\ell}(\xx)=\frac{1}{\sqrt{n_{\ell-1}}}W_{\ell}\xx_{\ell-1}(\xx),\quad \xx_{\ell-1}(\xx)=\sigma(g_{\ell-1}(\xx)),\quad \xx_0(x)=\xx
\]
for $\ell=1,2,\cdots,L$ where $n_{\ell}$ and $W^{\ell}\in\R^{n_{\ell}\times n_{\ell-1}}$ are the number of neurons and the weight matrix in $L$-the layer respectively. We suppose that all weights are initialized with standard Gaussian distributions independently. Let us further assume the circumstance that we aim to minimize a generic differentiable loss function $\mathcal{L}$:
\[
    \frac{{\rm d}\boldsymbol{\theta}(t)}{{\rm d}t}=-\nabla_{\tth}\mathcal{L}(\boldsymbol{\theta}(t))=-\sum_{i=1}^m\frac{\partial\mathcal{L}(y_i,z)}{\partial z}\bigg|_{z=f(\xx_i;\boldsymbol{\theta}(t))}\nabla_{\boldsymbol{\theta}}f(\xx_i;\boldsymbol{\theta}(t)),
\]
where $\boldsymbol{\theta}$ is a collection of all weights $\{W_{\ell}\}$. One of the main results of the paper \cite{NTK} is encapsulated in the following theorem.
\begin{theorem}
    Assume that the activation function $\sigma$ is of class $C^2$ and Lipschitz, and the loss function $\mathcal{L}$ is Lipschitz. Then the neural tangent kernel $\widehat{\boldsymbol{\Theta}}(t)$ converges in probability to  a deterministic kernel independent of $t>0$ as $n_1$, $n_2,\cdots n_{L-1}\rightarrow\infty.$
\end{theorem}
The above theorem demonstrates that training a neural network under Gaussian parameterization is equivalent to conducting a certain kernel method as widths go to
infinity.

\subsection{Neural Tangent Kernel for PINNs}\label{subsec:ACR}
In \cite{pre_work_1}, the authors investigated the NTK theory for PINNs. They have shown that when the parameters are initialized by the Gaussian distribution and the number of neurons in each layer goes to infinity, the NTK for PINNs is determined by the initial state and it remains constant during the training. To illustrate the relevant theory in more detail, let us consider the problem \eqref{main_eq}, and the corresponding physics-informed loss function \eqref{total_loss} with \eqref{res_loss}-\eqref{bdry_loss}. Then it was proved in \cite{pre_work_1} that the evolution of $u$ and $\mathcal{D}[u]$ can be described by the gradient flow:
\[
    \left[\begin{matrix}
        \frac{{\rm{d}}u(\boldsymbol{x}_b;\boldsymbol{\theta}(t))}{{\rm{d}}t}\\
        \frac{{\rm{d}}\mathcal{D}[u](\boldsymbol{x}_r;\boldsymbol{\theta}(t))}{{\rm{d}}t}
    \end{matrix}\right]=-\left[\begin{matrix}
        \boldsymbol{K}_{uu}(t)&\boldsymbol{K}_{ur}(t)\\
        \boldsymbol{K}_{ru}(t)&\boldsymbol{K}_{rr}(t)
    \end{matrix}\right]\left[\begin{matrix}
        u(\boldsymbol{x}_b;\boldsymbol{\theta}(t))-g(\boldsymbol{x}_b)\\
        \mathcal{D}[u](\boldsymbol{x}_r;\boldsymbol{\theta}(t))-f(\boldsymbol{x}_r)
    \end{matrix}\right],
\]
where $\boldsymbol{K}_{uu}(t)\in\mathbb{R}^{N_b\times N_b}$, $\boldsymbol{K}_{rr}\in\mathbb{R}^{N_r\times N_r}$ and $\boldsymbol{K}_{ru}(t)=[\boldsymbol{K}_{ur}(t)]^T\in \mathbb{R}^{N_r\times N_b}$ whose $(i,j)$-th entries are
\begin{equation}\label{NTK_PINN}
    \begin{aligned}
        &(\boldsymbol{K}_{uu})_{ij}(t) = \left<\frac{{\rm{d}}u(\boldsymbol{x}_b^i;\boldsymbol{\theta}(t))}{{\rm{d}}\boldsymbol{\theta}},\frac{{\rm{d}}u(\boldsymbol{x}_b^j;\boldsymbol{\theta}(t))}{{\rm{d}}\boldsymbol{\theta}}\right>,\\
        &(\boldsymbol{K}_{ru})_{ij}(t) = \left<\frac{{\rm{d}}\mathcal{D}[u](\boldsymbol{x}_r^i;\boldsymbol{\theta}(t))}{{\rm{d}}\boldsymbol{\theta}},\frac{{\rm{d}}u(\boldsymbol{x}_b^j;\boldsymbol{\theta}(t))}{{\rm{d}}\boldsymbol{\theta}}\right>,\\
        &(\boldsymbol{K}_{rr})_{ij}(t) = \left<\frac{{\rm{d}}\mathcal{D}[u](\boldsymbol{x}_r^i;\boldsymbol{\theta}(t))}{{\rm{d}}\boldsymbol{\theta}},\frac{{\rm{d}}\mathcal{D}[u](\boldsymbol{x}_r^j;\boldsymbol{\theta}(t))}{{\rm{d}}\boldsymbol{\theta}}\right>.
    \end{aligned}
\end{equation}
The matrix
\[\boldsymbol{K}(t)=\left[\begin{matrix}
        \boldsymbol{K}_{uu}(t)&\boldsymbol{K}_{ur}(t)\\
        \boldsymbol{K}_{ru}(t)&\boldsymbol{K}_{rr}(t)
    \end{matrix}\right]\]
is called the NTK of a PINN. Motivated by the fact that an infinitely-wide neural network becomes a Gaussian process and its NTK is constant during the training \cite{NTK}, \cite{pre_work_1} showed that, if the number of neurons of each layer goes to infinity, the NTK of a PINN $\boldsymbol{K}(t)$ indeed converges in probability to a deterministic kernel $\boldsymbol{K}^{*}$ and $\sup_{t\in[0,T]}\|\boldsymbol{K}(t)-\boldsymbol{K}(0)\|_2\rightarrow 0$ for any $T>0$, under some appropriate assumptions.

Since the kernel $\boldsymbol{K}(t)$ is positive semi-definite, by the spectral decomposition, we can see that the eigenvalues of the kernel $\boldsymbol{K}(t)$ provide some information on the convergence rate for training as presented in \eqref{2_2_4} for given data points (see \cite{pre_work_1} for details). However, as the number of training data increases, it is not easy to explicitly compute all the eigenvalues individually. Therefore, based on the fact that the trace is equal to the sum of eigenvalues, motivated by \cite{pre_work_1}, we shall used the following averaged quantity to measure the approximate rate of convergence for the training of PINNs
\begin{equation}\label{avg_conv_rate}
\frac{{\rm Tr}({\boldsymbol{K}_{uu}}(t))}{N_b}+\frac{{\rm Tr}({\boldsymbol{K}_{rr}}(t))}{N_r}.
\end{equation}
In the later section, we will explicitly compute the above quantity when the variables are scaled to illustrate the faster convergence of the proposed method in this paper.


\section{Variable-Scaling Physics-Informed Neural Networks}\label{sec:method}
In this section, we will describe the method we mainly propose in this paper.  
The key idea of our method is to use a variable-scaling technique in the training process. More precisely, we exploit the change of variable $\xx\mapsto\overline{\xx}/N$ to define a loss function, in order to relieve the potential stiff behavior of solutions. After the training is completed, we scale back to the original variable for the solution prediction.

\begin{figure}[H]
\begin{center}
\includegraphics[width=0.75\textwidth]{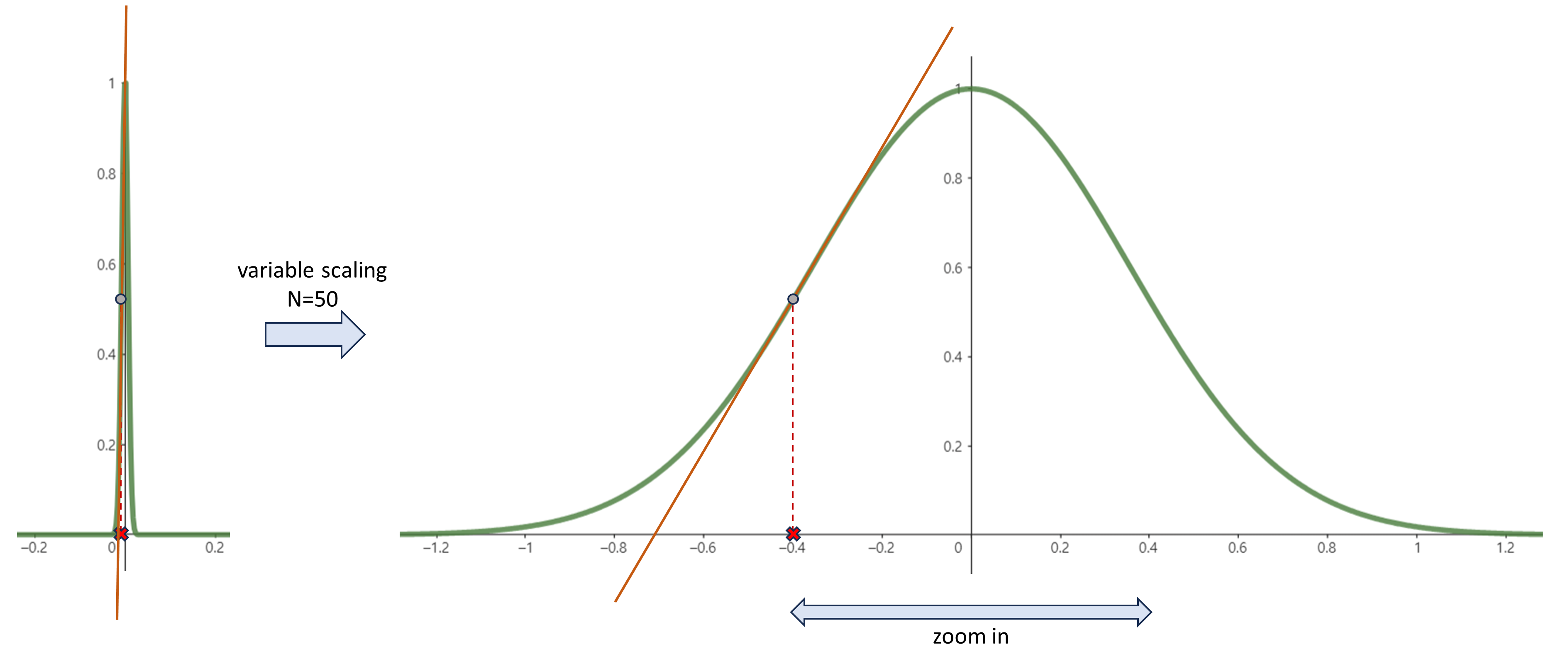}
\end{center}
\caption{An illustration of the effect of the variable scaling when the scaling factor $N=50$ is applied to the function $\exp( -10000x^2)$. The scaling has the effect of ``zoom-in'' and the resulting function exhibits relatively less stiff behavior on the designated domain.}
\label{VS_ex}
\end{figure}

To illustrate the methodology in more detail, let us consider the following PDE of the form:
\begin{equation}\label{given_PDE}
\begin{aligned}
    F(\xx;u,Du,D^2u,\cdots,D^ku)&=0 &&\xx\in\Omega,\\
    u(\xx)&=g(\xx) && \xx\in\partial\Omega,
\end{aligned}
\end{equation}
where $F(\xx)=\mathcal{D}[u](\xx)-f(\xx)$ in \eqref{main_eq} and $D^j$ denotes the mixed partial derivatives of degree $j\in\mathbb{N}$. Now the corresponding residual loss function for the PINN is then given by
\begin{equation}\label{PINN_loss}
    \mathcal{L}_{\rm original}=\frac{\lambda_r}{N_r}\sum^{N_r}_{i=1}\bigg|F(\xx^i_r;u(\xx^i_r),Du(\xx^i_r),D^2u(\xx^i_r),\cdots,D^ku(\xx^i_r))\bigg|^2+\frac{\lambda_b}{N_b}\sum^{N_b}_{j=1}\left|u(\xx^j_b)-g(\xx^j_b)\right|^2,
\end{equation}
where the $\xx^i_r$ and $\xx^j_b$ are i.i.d. random samples chosen from the uniform distribution $\xx^i_r\sim\mathcal{U}(\Omega)$ and $\xx^j_b\sim\mathcal{U}(\partial\Omega)$ respectively. Next, let us introduce a new variable $\overline{\xx}=N\xx\in N\Omega:=\{N\xx:\xx\in\Omega\}$ with the notation $u(\overline{\xx}/N)=v(\overline{\xx})$. Intuitively, as we can see from Figure \ref{VS_ex}, the change of variable $\xx=\overline{\xx}/N$ can be interpreted as an $N$-times ``zoom-in'' for the solution profile, and thus it has an effect of reducing any sharp transitions of solution dynamics on the designated domain. 
By the chain rule, we know that
\[
D_{\xx}u(\xx)=D_{\overline{\xx}}u\left(\frac{\overline{\xx}}{N}\right)\frac{{\rm d}\overline{\xx}}{{\rm d}\xx}=ND_{\overline{\xx}}v(\overline{\xx}),
\]
and this can be easily extended to $D^k_{\overline{\xx}}u(\xx)=N^kD_{\overline{\xx}}v(\overline{\xx})$. Therefore, the given PDE \eqref{given_PDE} can be rewritten as
\begin{equation}\label{rewritten_PDE}
\begin{aligned}
    F(\overline{\xx}/N;v,NDv,N^2D^2v,\cdots,N^kD^ku)&=0 &&\overline{\xx}\in N\Omega,\\
    v(\overline{\xx})&=g(\overline{\xx}/N) && \xx\in N\partial\Omega.
\end{aligned}
\end{equation}
The associate loss function for the representation \eqref{rewritten_PDE} is
\begin{equation}\label{scaled_loss}
    \mathcal{L}_{\rm scaled}=\frac{\lambda_r}{N_r}\sum^{N_r}_{i=1}\bigg|F(\overline{\xx}^i_r/N;v(\overline{\xx}^i_r),NDv(\overline{\xx}^i_r),N^2D^2v(\overline{\xx}^i_r),\cdots,N^kD^kv(\overline{\xx}^i_r))\bigg|^2+\frac{\lambda_b}{N_b}\sum^{N_b}_{j=1}\bigg|v(\overline{\xx}^j_b)-g\bigg(\frac{\overline{\xx}^j_b}{N}\bigg)\bigg|^2,
\end{equation}
with $\overline{\xx}^i_r\sim\mathcal{U}(N\Omega)$ and $\overline{\xx}^j_b\sim\mathcal{U}(N\partial\Omega)$. In this scaled loss function, the weight parameters $\lambda_r$ and $\lambda_b$ should be carefully chosen according to the result of the scaling. For example, they can be chosen to cancel a common factor or extreme coefficients (see Section \ref{sec:exp} for details). The key novel part of our method is to minimize the loss function \eqref{scaled_loss} instead of \eqref{PINN_loss}. Once the training is finished and we obtain the solution prediction $v(\overline{\xx})$, we pull the variable back to the original scale by the formula
\begin{equation}\label{zoomout}
    u(\xx)=v(N\xx),
\end{equation}
which is the original solution we aimed to predict at the initial stage. Since the domain $\Omega$ has become larger to $N\Omega$ in the method, one may wonder whether we need more random samples for the training. However, as we will illustrate in the numerical experiments in the next section, this is not the case, and we can significantly improve the training efficiency even if we train the model with the same number of random samples from the scaled domain.

\section{Numerical Experiments}\label{sec:exp}

In this section, we will perform a series of numerical experiments to evaluate the performance of the proposed training method against the standard PINNs. As mentioned earlier, the choice of parameters $\lambda_r$ and $\lambda_b$ in \eqref{scaled_loss} as well as the scaling factor $N$ is important for the overall performance and they should be carefully selected. As we will see in the numerial experiments in this section, as $N$ increases, the parameters related to training become more sensitive, which tends to make the training somewhat unstable with more fluctuation in learning curves. Therefore, to leverage the advantages of the VS-PINN while ensuring adequate training stability, it is necessary to choose an appropriate $N$. We will also discuss this issue in more details in the later section of this paper, and here we shall demonstrate some basic strategies for choosing these parameters in the experiments.

 Throughout the section, the model trained with the original PINN will be compared with the model trained by the proposed variable-scaling technique, and all the other settings including the model architecture, initialization, and the optimizer will be the same.  
And most importantly, even though the domain is $N$-times extended, the number of collocation points for the training is the same for both the original method and the variable-scaling method. 

\subsection{Wave Equation}

As a starting point, we shall conduct some experiments solving the one-dimensional wave equation, where the initial data with a high-frequency part is imposed. Namely, we will solve the following problem:

\begin{equation}\label{wave_problem}
\begin{aligned}
&u_{tt}-u_{xx}=0 && (x,t)\in(0,1)\times (0,1),\\
&u(0,t)=u(1,t)=0&&t\in[0,1],\\
&u(x,0)=0&&x\in[0,1],\\
&u_t(x,0)= 2\pi\sin(2\pi x)+10\pi\sin(10\pi x)&&x\in[0,1].
\end{aligned}
\end{equation}
The exact solution for the above problem can be represented as
\[u(x,t) = \sin(2\pi x)\sin(2\pi t)+\sin(10\pi x)\sin(10\pi t),\] 
where the $\sin(10\pi x)\sin(10\pi t)$ part exhibit a high frequency.

\begin{figure}
    \centering
  \includegraphics[width=0.4\linewidth]{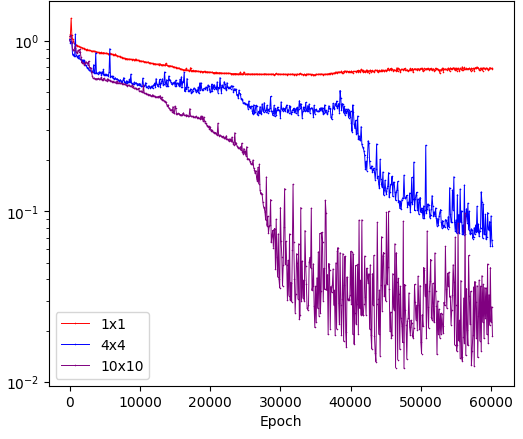}
    \caption{Learning curves for the VS-PINNs with the scales $N=1,4,10$.}
    \label{wave_learn}
\end{figure}

For the baseline model, we use a neural network with 4 hidden layers with 128 neurons for each, activated by the hyperbolic tangent function. By scaling the variables to expand the domain, we can relatively mitigate the high-frequency regions, which could potentially improve the training process. To demonstrate this, we conducted experiments with three different scaling factors $N=1,4,10$, which leads us to compute the solutions in three different computational domains $[0,1]^2,[0,4]^2$ and $[0,10]^2$. The scaled problem of \eqref{wave_problem} can be written as
\begin{equation}\label{wave_eq}
\begin{aligned}
&N^2u_{tt}-N^2u_{xx}=0 &&(x,t)\in(0,N)\times (0,N),\\
&u(0,t)=u(N,t)=0&&t\in[0,N],\\
&u(x,0)=0&&x\in[0,N],\\
&u_t(x,0)= \frac{2\pi}{N}\sin\left(\frac{2\pi x}{N}\right)+\frac{10\pi}{N}\sin\left(\frac{10\pi x}{N}\right)&&x\in[0,N],
\end{aligned}
\end{equation}
and the corresponding scaled loss is
\begin{equation}\label{wave_loss}
    \mathcal{L}_{\rm{scaled}} = 2\mathcal{L}_{\rm data}+\frac{1}{N^4}\mathcal{L}_{\rm res}.
\end{equation}
As mentioned earlier, note that it is reasonable to choose  $\lambda_r=\frac{1}{N^4}$ to cancel the $N^4$ factor in the residual loss function. For the training, we randomly sample 6400 collocation points for $\mathcal{L}_{\rm res}$ and 512 points for $\mathcal{L}_{\rm data}$. The learning-curve comparison for the methods with different scales $N=1$ (standard PINN) and $N=4,10$ (VS-PINN) is presented in Figure \ref{wave_learn}. For each epoch, we compute the relative $L^2$ error between the true solution and the PINN to evaluate the training efficiency. As the scale increases, we can observe in the figure that the error decreases more quickly and rapidly, meaning that the VS-PINN can learn the target solution with high frequency more efficiently.

\begin{figure}[H]
\centering
\includegraphics[width=0.9\linewidth]
{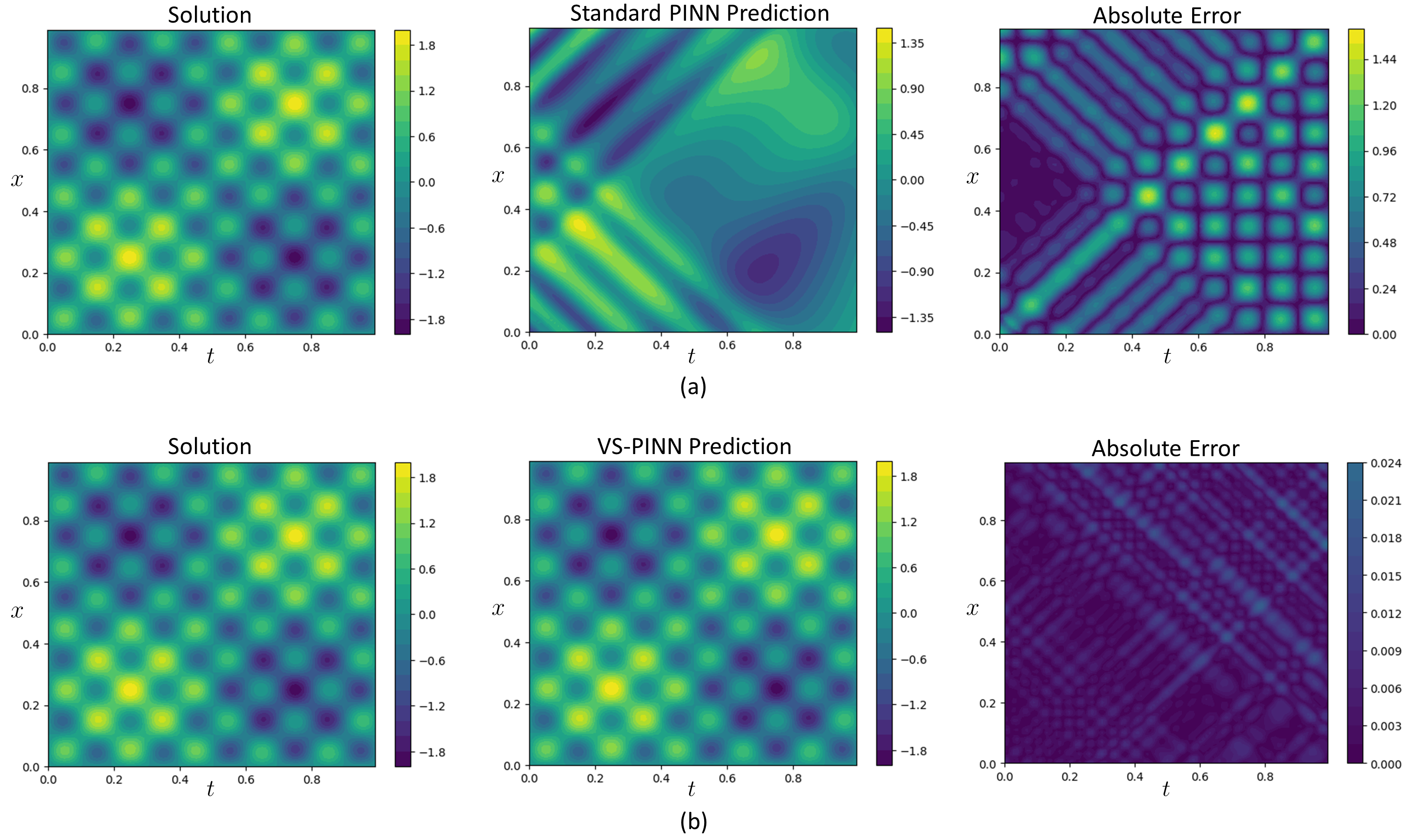}
\caption{One-dimensional wave equation: (a) The exact solution versus the predicted solution by training the standard PINN. The relative $L^2$ error is 6.31e-1. (b) The exact solution versus the predicted solution by training the VS-PINN with the scale $N=10$. The relative $L^2$ error is 1.20e-2.}
\label{wave_pred}
\end{figure}

We train the model in a zoom-in state, followed by an evaluation in a zoom-out state with the process: $u(x_{\rm test},t_{\rm test})=u_{\rm model}(Nx_{\rm test},Nt_{\rm test})$.
After the 60000 epochs of training, a comparison between the true solution and the predicted solution is illustrated in Figure \ref{wave_pred}. As we can see from Figure \ref{wave_pred}, the solution prediction with variable scaling is much more accurate for the high-frequency solution compared to the solution computed by the standard PINN.

\subsection{Allen--Cahn Equation}
We shall also apply the variable-scaling technique to the PINNs for other types of equations. In particular, it is noteworthy that the change of variable can be used to offset the effect of the extreme coefficient that appeared in the given equation. This idea will be examined through the numerical experiments for the Allen--Cahn equation in this subsection and the singular perturbation problem in the next subsection in order.

The Allen–Cahn equation arises in mathematical physics, which describes the process of phase separation in multi-component alloy systems. From a mathematical perspective, it is the reaction-diffusion-type equation with nonlinear terms, typically assumed to have a small diffusion parameter. We shall perform the experiment in the same setting adopted in \cite{pinn01}, which is represented by the following formulation:
\begin{equation}
    \begin{aligned}
        &u_t-0.0001u_{xx}+5u^3-5u = 0 &&x\in(-1,1),\;\;\;t\in(0,1),\\
        &u(x,0)=x^2\cos(\pi x) && x\in[-1,1],\\
        &u(-1,t)=u(1,t)  && t\in[0,1],\\
        &u_x(-1,t)=u_x(1,t)  && t\in[0,1].
    \end{aligned}
\end{equation}
For the model architecture, we exploit the neural network consisting of 4 hidden layers with 64 neurons for each, and the corresponding parameters are initialized with random initialization using a fixed seed. For training, we first proceed with supervised learning and train the neural network with pre-computed training data as performed in \cite{pinn01}. In addition, in order to further highlight the effectiveness of our method, we will also train the VS-PINN in an unsupervised manner and confirm that we can obtain satisfactory results even for unsupervised learning.

For the supervised learning part, data points that are randomly chosen from the exact solution at $t = 0.1$ as conducted in \cite{pinn01}. For the other collocation points (both for supervised and unsupervised learning) we select $128$ random samples for the periodic boundary conditions, $200$ points for the initial condition and $10000$ data points for the interior.
For the scaling, we consider the change of variable only for the spatial variable $x\mapsto\overline{x}/N$. Through a derivative, we can cancel the diffusion coefficient which might cause a stiff behavior. More precisely, to define the scaled loss function, we shall consider the following scaled problem:
\begin{equation}\label{AC_scaling}
    \begin{aligned}
        &u_t-N^2\cdot0.0001u_{xx}+5u^3-5u = 0 &&x\in(-N,N),\;\;\;t\in(0,1),\\
        &u(x,0)=(x/N)^2\cos(\pi x/N) &&x\in[-N,N],\\
        &u(-N,t)=u(N,t)&& t\in[0,1],\\
        &u_x(-N,t)=u_x(N,t)&&t\in[0,1].
    \end{aligned}
\end{equation}
As we can see from the scaled problem \eqref{AC_scaling}, it would be reasonable to choose some large scaling factor $N$ so that we can cancel the diffusion coefficient $0.0001$, which mitigates the potential extreme behavior of solution caused by the biased coefficient. In the simulation, we shall examine the cases when $N$ is $1$ (the standard PINN), $10$ and $100$ (the VS-PINN). Finally, we use the training loss function with adjusted parameters $\lambda_{\rm data}=2$ and $\lambda_{\rm res}=0.3$:
\[\mathcal{L}_{\rm scaled}=2\mathcal{L}_{\rm data}+0.3\mathcal{L}_{\rm res}.\]
As before, the model is trained in a zoom-in state and the performance is evaluated in a zoom-out state which is the original scale: $u(x_{\rm test},t_{\rm test})=u_{\rm model}(Nx_{\rm test},t_{\rm test})$. 

\begin{figure}
    \centering
    \begin{subfigure}{0.49\textwidth}
        \centering
        \includegraphics[width=\linewidth]{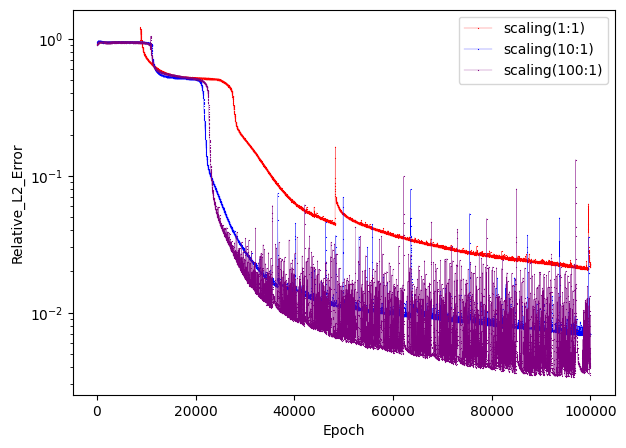}
        \caption{Supervised learning}
        \label{AC_super}
    \end{subfigure}    
    \begin{subfigure}{0.49\textwidth}
        \centering
        \includegraphics[width=\linewidth]{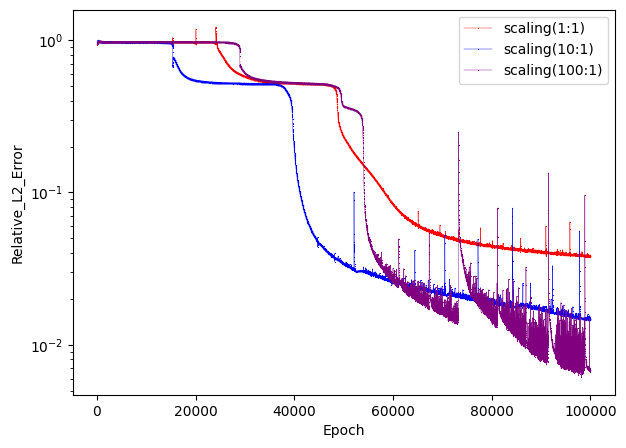} 
        \caption{Unsupervised learning}
        \label{AC_unsuper}
    \end{subfigure}
\caption{Traning dynamics for the Allen--Cahn equation using VS-PINN with the scales $N=1$ (standard PINN), $10$, $100$. (a) is for the model trained with the pre-computed training data at $t=0.1$ and (b) is for the model trained in an unsupervised manner.}
\label{AC_LC}
\end{figure}

\begin{figure}[H]
    \centering
    \includegraphics[width=1\textwidth]{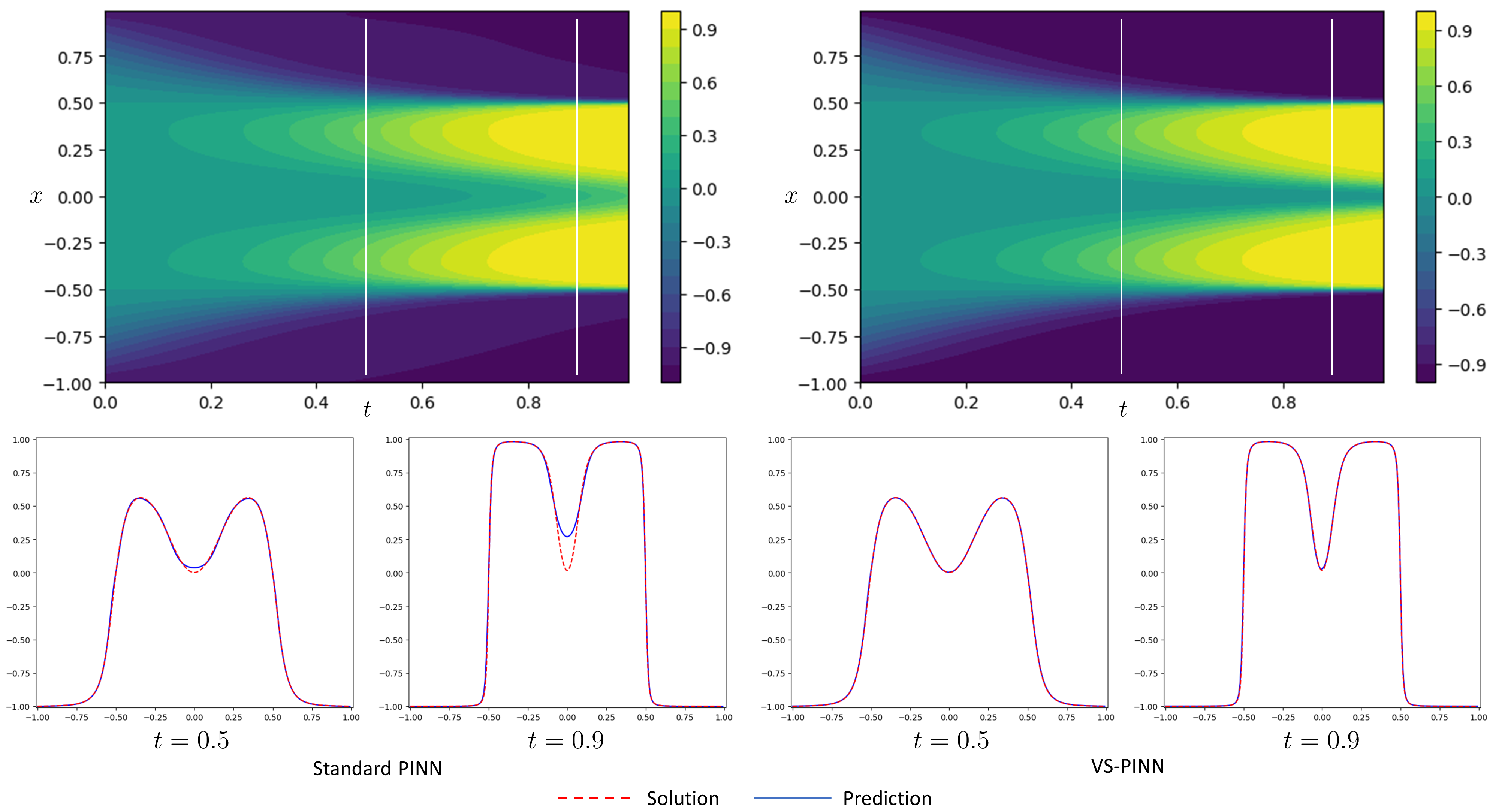}
    \caption{{\it{Top left}}: Plot of the prediction $\widehat{u}(x, t)$ by the standard PINN. {\it{Top right}}: Plot of the prediction $\widehat{u}(x, t)$ via the VS-PINN. {\it{Bottom left}}: Snapshots of the prediction $\widehat{u}(x, t)$ using the original PINN (blue solid curve) vs. the exact solution $u(x, t)$ (red dashed curve) at $t=0.5$ and $t=0.9$. {\it{Bottom right}}: Snapshots of the prediction $\widehat{u}(x, t)$ via the VS-PINN (blue solid curve) vs. the exact solution $u(x, t)$ (red dashed curve) at $t=0.5$ and $t=0.9$.}
    \label{AC_pred}
\end{figure}

Figure \ref{AC_LC} displays the relative $L^2$ error for the standard PINN and the VS-PINN as a function of the training epoch. From the training curves, we can confirm that we can significantly improve the training efficiency and obtain higher accuracy both for supervised and unsupervised scenarios. In the case of unsupervised learning, as presented in Figure \ref{AC_LC} (b), it was observed that as the domain expanded, the point at which the error began to decrease was delayed when $N=100$. However, as the training proceeds further, we can finally obtain a faster convergence with higher accuracy.  For both cases, the solution profiles trained in an unsupervised manner are presented in Figure \ref{AC_pred}. In particular, the snapshots at $t=0.5$ and $t=0.9$ are illustrated in the bottom line for both cases. The error was observed to be 2.05e-2 for the standard PINN ($N=1$) and 3.36e-3 for the VS-PINN ($N=100$). As we can see from these pictures, the solution prediction via the VS-PINN better captures the stiff behavior of the solution. 

\subsection{Boundary-Layer Problem}
The next example focuses on the singularly perturbed boundary-layer problem, which is a well-known challenge in scientific computing. If the diffusion coefficient is small, the solution exhibits a steep increase in thin regions near boundaries which may lead to failures in overall prediction. Typically, in such a case, neural-network-based solution prediction often has difficulties in finding optimal solutions. This is also mainly due to the spectral bias phenomenon mentioned in the introduction of this paper \cite{bias_1}. 
We aimed to alleviate such sharp changes through a variable scaling. Specifically, as mentioned earlier, if we make a change of variable, a small diffusion coefficient can be canceled by the Jacobian factor. To illustrate this idea, we shall consider the following one-dimensional singular perturbation problem: for a small parameter $0<\varepsilon\ll1$,

\begin{equation}\label{sp_prob}
\begin{aligned}
&\varepsilon u_{xx}+u_x+1=0\;\;\;x\in[0,1],\\
&u(0)=u(1)=0.
\end{aligned}
\end{equation}
The exact solution can be computed analytically, and is presented as follows:
\[
u(x) = \frac{1}{1-e^{-1/\varepsilon}}-x-\frac{1}{1-e^{-1/\varepsilon}}e^{-x/\varepsilon}.
\]
When $\varepsilon$ is sufficiently small, we can ignore the term $e^{-1/\varepsilon}$ and the solution is reduced to
\[u(x) = 1-x-e^{-x/\varepsilon}.\]
This solution has a slope of $1/\varepsilon-1$ at $x=0$ and it attains its maximum at $x=-\varepsilon\ln(\varepsilon)$, indicating there is an explosive change near the boundary.
We set $\varepsilon=$1e-6 and perform the change of variable with the scale $N=1000$ so that we can remove the effect of $\varepsilon>0$ by multiplying the Jacobian factor. Then, the original equation is transformed as follows:
\begin{equation}\label{SP}
    \begin{aligned}
        &u_{xx}+1000u_x+1=0\;\;\;x\in[0,1000],\\
        &u(0)=u(1000)=0,
    \end{aligned}
\end{equation}
and we set the loss as
\[\mathcal{L}_{\rm scaled}=20\mathcal{L}_{\rm data}+\mathcal{L}_{\rm res}.\]

\begin{figure}
    \centering
    \includegraphics[width=0.8\textwidth]{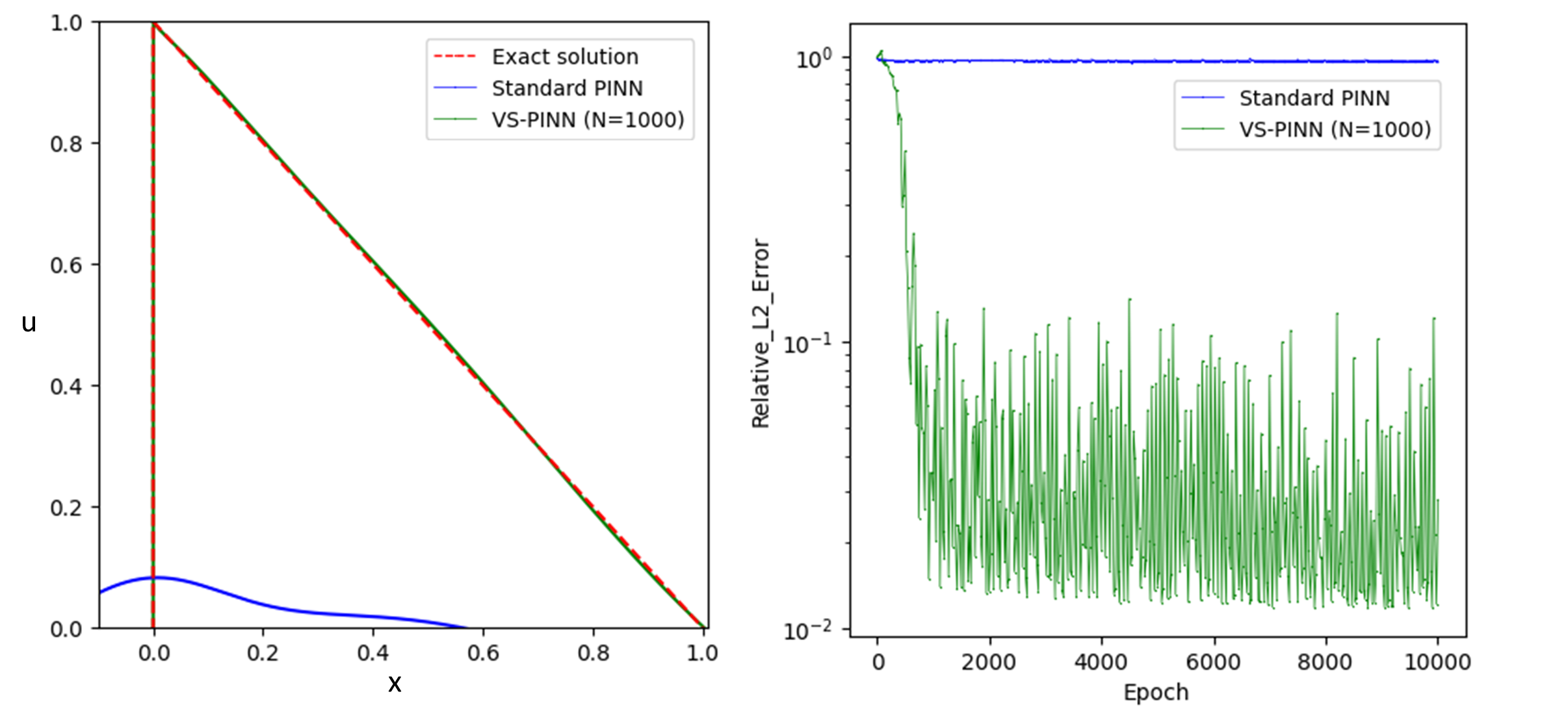}
    \caption{{\it{Left}}: Exact solution of \eqref{sp_prob} compared with the solution prediction by the standard PINN and the VS-PINN. {\it{Right}}: Learning curves for both the original PINN and the VS-PINN.}
\end{figure}

We shall use the model with 8 hidden layers, with 20 neurons per layer and the model is initialized with a fixed seed and is trained using 1000 interior collocation data. When the model is trained with a variable scaling, the highest accuracy achieved was 1.17e-2, while in the original case, there is almost no progress in training and the learning curve remains mostly unchanged. Therefore, we can conclude that the mitigation of stiff behavior through the variable scaling has a significant impact on training and can make a much more accurate solution prediction for the boundary-layer problem.

\subsection{The Navier--Stokes Equations}
Another advantage of using the VS-PINN is that the neural network can learn the nonlinear structure more efficiently by alleviating the potentially complicated dynamics of solutions. As a last experiment, we will discuss the experimental results of the challenging Navier--Stokes equations. The numerical computation of solutions to fluid equations has a long and rich history, and it has seen significant advancements across various fields. Furthermore, in recent years, there has been a wide range of studies to predict fluid flow based on machine-learning methods, and this topic has become of independent interest. For example \cite{NS_comp} proposed a mixed-variable PINN to compute the solution of the incompressible Navier--Stokes equations, \cite{pinn_fluid_4} proposed the enhanced PINN to solve the three-dimensional Navier--Stokes equations and \cite{pinn_fluid_3} performed PINN computations by considering both the velocity-pressure formulation and the vorticity-velocity to simulate the flows from laminar to turbulent formulation.

\begin{figure}
    \centering    \includegraphics[width=0.65\textwidth]{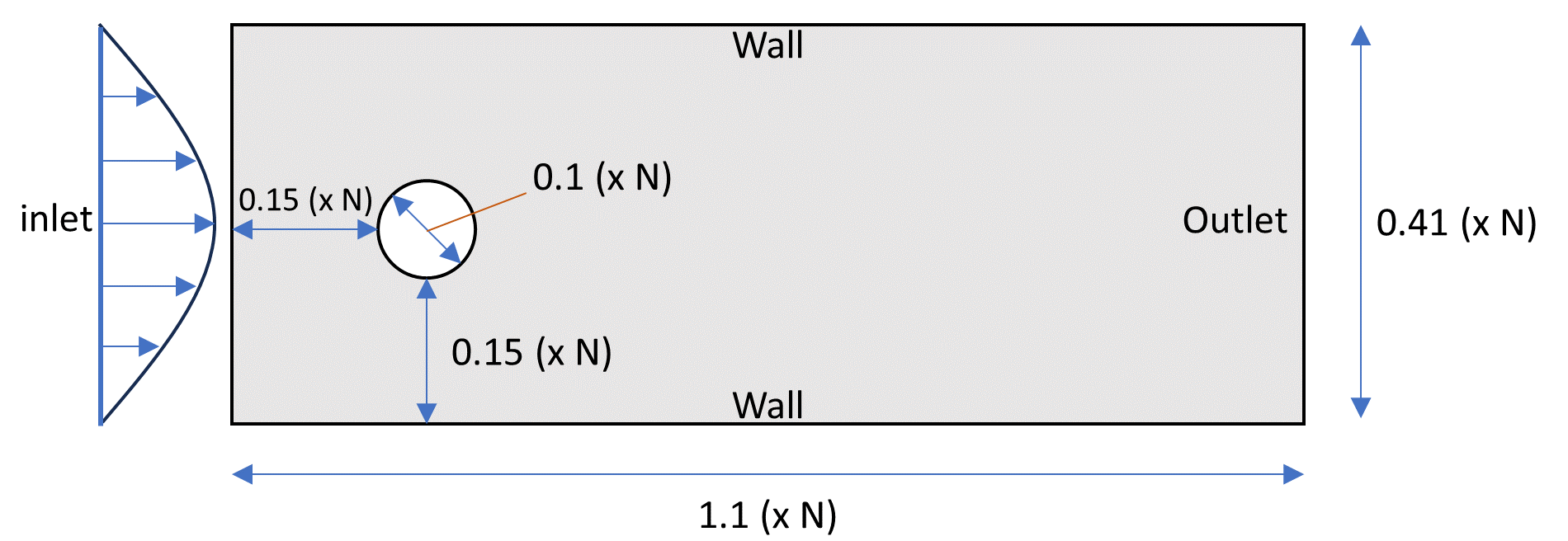}
    \caption{An illustration of the computational domain where we compute the solution of the incompressible Navier--Stokes equations.}
    \label{fig:enter-label}
\end{figure}

As an extension of these prior research papers, in this paper, we shall conduct the experiments for the steady and transient laminar flows passing over a circular cylinder in the two-dimensional domain without any training data as investigated in \cite{NS_comp}. 
More precisely, we shall consider the following incompressible Newtonian fluid flow problem, described by the following Navier--Stokes equations:\begin{equation}\label{NS_equation_ori}
    \begin{aligned}
        u_x + v_y &= 0 &&(x,y)\in\Omega, \\
        uu_x+vu_y &= -\frac{1}{\rho}p_x + \frac{\mu}{\rho}(u_{xx}+u_{yy})&&(x,y)\in\Omega,\\
        uv_x+vv_y &= -\frac{1}{\rho}p_y + \frac{\mu}{\rho}(v_{xx}+v_{yy})&&(x,y)\in\Omega,
    \end{aligned}
\end{equation}
where $(u,v)\in\R^2$ is the velocity vector, $p\in\R$ is the pressure and the viscosity and the density of the fluid are set to be $\mu=0.02$ and $\rho=1$ respectively. The computational domain $\Omega$ is depicted in Figure \ref{fig:enter-label} which will be expanded by a factor of $N$ in the training for the VS-PINN.

A velocity profile is applied on the inlet while the condition of zero pressure is imposed on the outlet, as shown in Figure \ref{fig:enter-label}. No-slip boundary conditions are assumed on the wall and cylinder boundaries. Namely, we enforce the following conditions to the equations \eqref{NS_equation_ori}:

\begin{equation}\label{NS_bdry}
    \begin{aligned}
        &u(x,y) = v(x,y) = 0 &&\text{on wall and cylinder boundaries}, \\
        &p(1.1,y) = 0 &&y\in[0,0.41],\\
        &u(0,y) = 4(0.41-y)\frac{y}{0.41^2} &&y\in[0,0.41], \\
        &v(0,y) = 0 &&y\in[0,0.41].
    \end{aligned}
\end{equation}

As before, the corresponding scaled boundary-value problem can be formulated as
\begin{equation}\label{NS_equation_sc}
    \begin{aligned}
        Nu_x + Nv_y &= 0 &&(x,y)\in N\Omega, \\
        uNu_x+vNu_y &= -\frac{1}{\rho}Np_x + \frac{\mu}{\rho}(N^2u_{xx}+N^2u_{yy}) &&(x,y)\in N\Omega,\\
        uNv_x+vNv_y &= -\frac{1}{\rho}Np_y + \frac{\mu}{\rho}(N^2v_{xx}+N^2v_{yy}) &&(x,y)\in N\Omega,\\
        u(x,y) &= v(x,y) = 0 &&\text{on wall and cylinder boundaries,} \\
        p(N\cdot1.1,y) &= 0 &&y\in[0,N\times0.41],\\
        u(0,y) &= 4(0.41-y/N)\frac{y/N}{0.41^2}&&y\in[0,N\times0.41], \\
        v(0,y) &= 0&&y\in[0,N\times0.41],
    \end{aligned}
\end{equation}
and we define the associated scaled loss function by
\[
    \mathcal{L} = \frac{1}{N^2}\mathcal{L}_{\rm res}+2\mathcal{L}_{\rm data}.
\]

\begin{figure}
    \centering
    \includegraphics[width=0.8\textwidth]{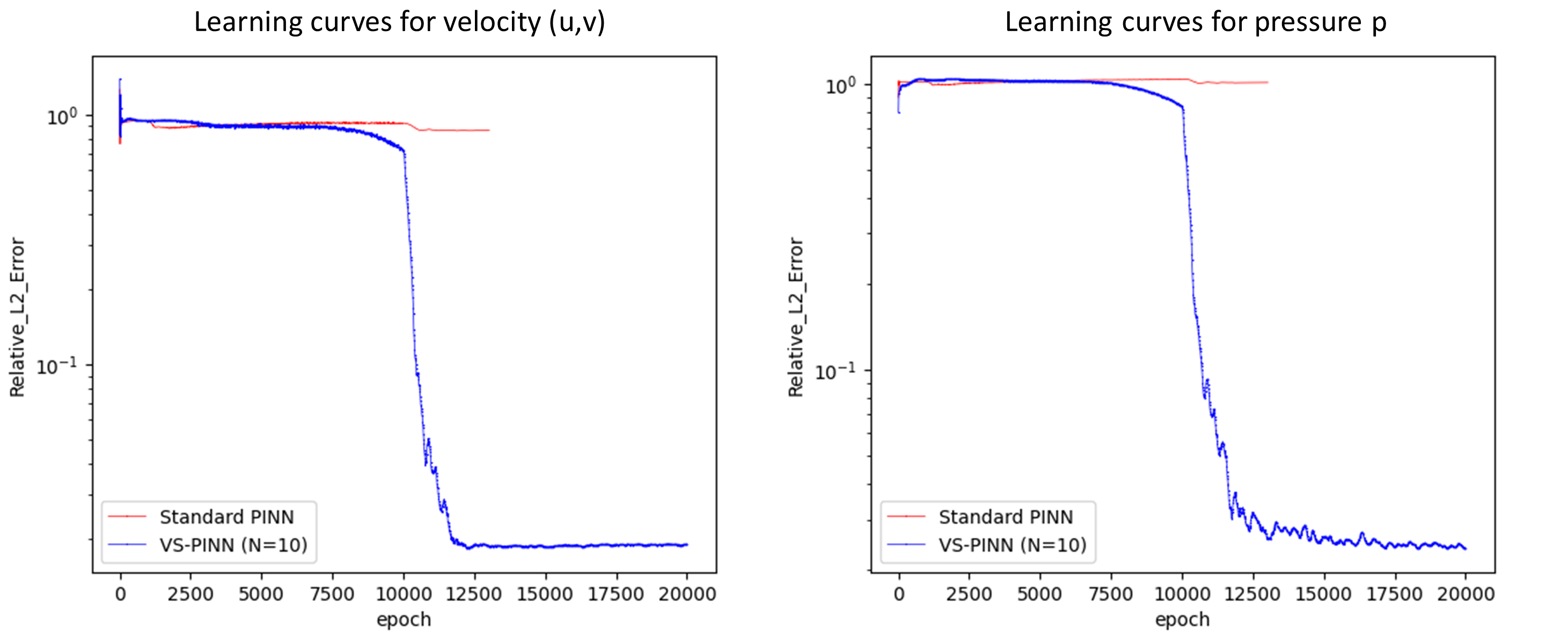}
    \caption{The standard PINN vs. the VS-PINN trained by the Adam+L-BFGS optimizer.}
    \label{NS_comb}
\end{figure}
We train the models both for the standard PINN and the VS-PINN with two different optimizing strategies to highlight the efficiency of the proposed training method. For the comparison, we first exploit the strategy used in \cite{NS_comp}, that is, Adam optimizer with first $10000$ epochs and followed by the extra $10000$ epochs using the L-BFGS. For the Adam optimizer, we randomly sampled $1200$ collocation points from the boundary and $6000 + 600$ (near the cylinder) interior data points. For the L-BFGS optimizer, 1200 random samples were chosen from the boundary and $8000 + 1000$ (near the cylinder) data points were used for the interior. The Adam optimizer randomly samples the points from the interior at each epoch, while L-BFGS samples are fixed throughout training. The test errors are 1.83e-2 for the velocity and 2.37e-2 for the pressure, and the training curves are presented in Figure \ref{NS_comb}. 
\begin{figure}
    \centering
    \includegraphics[width=0.8\textwidth]{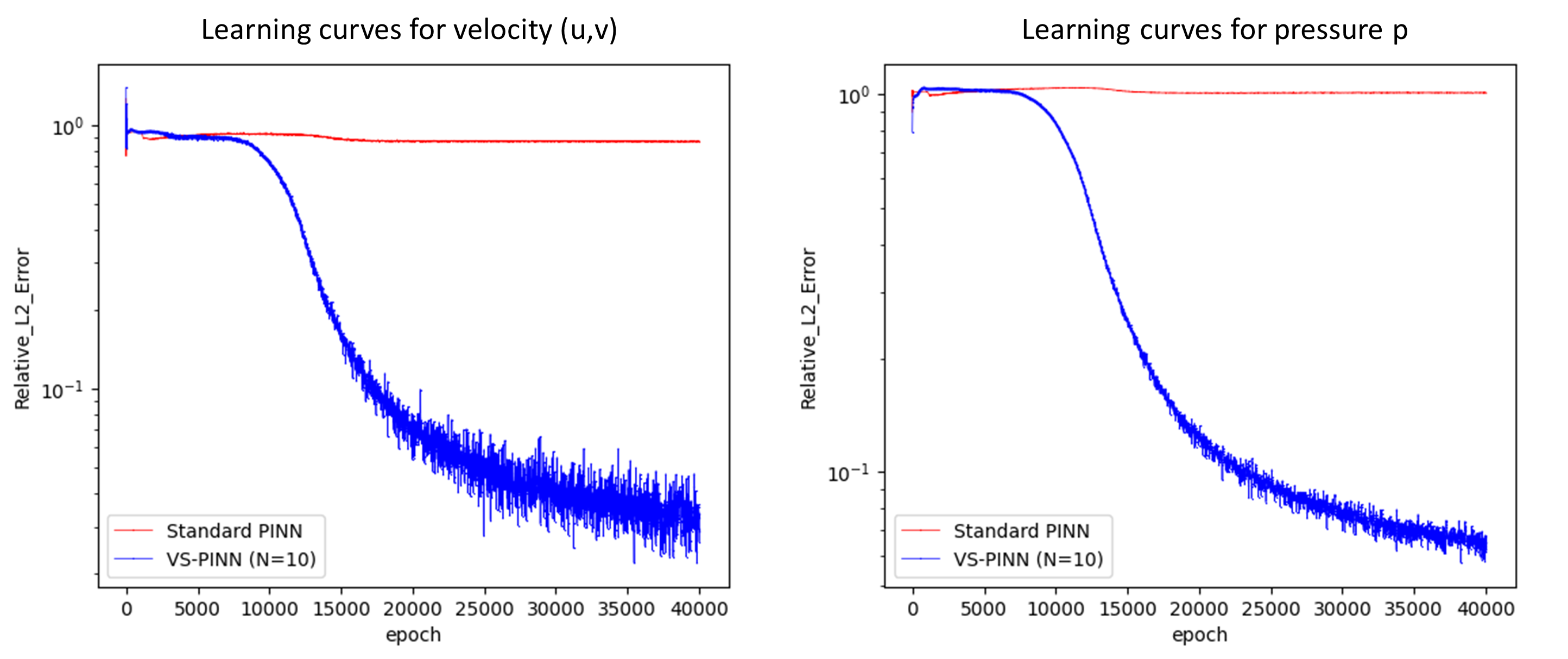}
    \caption{The standard PINN vs. the VS-PINN trained by the Adam optimizer only with a smaller number of collocation points.}
    \label{NS_adam}
\end{figure}

Note that even though we utilized the smallest model (5 hidden layers with 40 neurons for each) among the ones used in \cite{NS_comp}, the achieved accuracy is comparable to the accuracy obtained by using the larger model (7 hidden layers with 40 neurons each) in \cite{NS_comp}, which trained on $50000$ data points. Hence, we can confirm that the VS-PINN is effective in computing the approximate solutions for the Navier--Stokes equations.

On the other hand, the second strategy only utilizes the Adam optimizer with $1200$ collocation points on the boundary and $6000 + 600$ (near the cylinder) data points in the interior and proceeds $40000$ iteration. The final test errors are 2.18e-2 for the velocity and 5.75e-2 for the pressure. It is noteworthy that even though we only use much smaller number of training samples compared with the previous option using the L-BFGS as done in \cite{NS_comp}, we can still obtain the result with high accuracy, which was not observed in other previous works.

\begin{figure}
    \centering
    \includegraphics[width=1\textwidth]{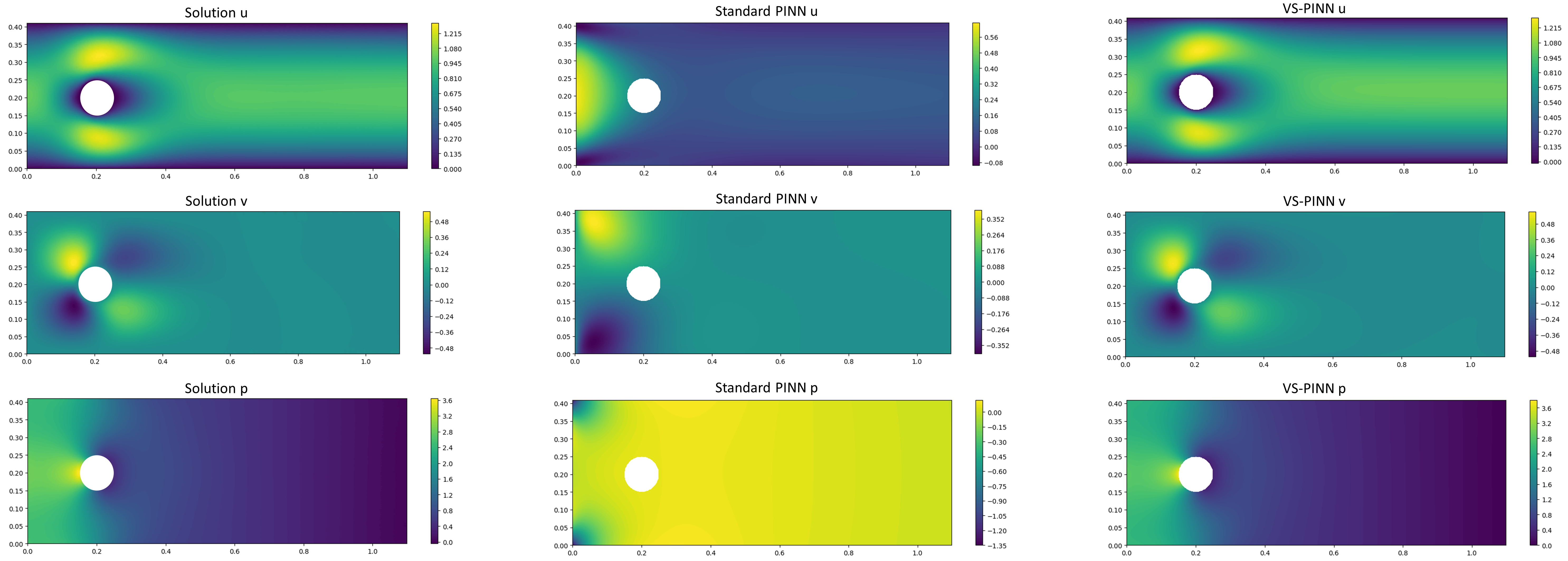}
    \caption{Contour plots for the solution and the approximate solutions. {\textit{Left column}}: Exact solutions for the problem \eqref{NS_equation_ori}-\eqref{NS_bdry}. {\textit{Middle column}}: Solution predictions using the standard PINN ($N=1$). {\textit{Right column}}: Solution predictions using the VS-PINN ($N=10$).}
    \label{NS_visual}
\end{figure}

The exact solutions and the predicted velocity and pressure computed by both the standard PINN and the VS-PINN are presented in Figure \ref{NS_visual}. From these contour plots, we can observe that the VS-PINN makes a much more accurate solution compared with the solution predicted by the original PINN. Moreover, since our method is universal, we can combine the proposed VS-PINN method with the previous ones to obtain better performance. For example, we can use various formulations proposed in the previous works \cite{pinn_fluid_3, pinn_fluid_4, NS_comp} and apply the variable-scaling technique to these additional formulations.

\section{Neural Tangent Kernel Analysis}\label{sec:thoery}

In this section, we shall investigate the effect of the variable scaling on training from an NTK perspective. Let us first describe the relevant problem that we shall consider for the analysis. As we discussed in  Section \ref{subsec:ACR}, the averaged convergence rate can be expressed in terms of the quantities ${\rm Tr}({\boldsymbol{K}_{uu}(0)})$ and ${\rm Tr}({\boldsymbol{K}_{rr}(0)})$. In order to convey the idea clearly, we will focus on a simple computational setting and explicitly compute the values ${\rm Tr}({\boldsymbol{K}_{uu}(0)})$ and ${\rm Tr}({\boldsymbol{K}_{rr}(0)})$. More precisely, we will consider the one-dimensional Poisson problem with a homogeneous Dirichlet boundary condition
\begin{equation}\label{theory_ori}
    \begin{aligned}
        -u_{xx} &= f(x) &&x\in [0,1],\\
        u(0)&=u(1)=0,
    \end{aligned}
\end{equation}
whose scaled problem can be written as
\begin{equation}\label{theory_scaled}
    \begin{aligned}
        -u_{xx} &= f(x/N)/{N^2}\;\; &&x\in [0,N],\\
        u(0)&=u(N)=0.
    \end{aligned}
\end{equation}
We then define the  corresponding scaled loss as
\[\mathcal{L}_{\rm scaled}=\frac{1}{N_r}\sum_{i=1}^{N_r}|u_{xx}(x_r^i;\boldsymbol{\theta})+f(x_r^i/N)/{N^2}|^2 + \frac{1}{N_b}\sum_{j=1}^{N_b}|u(x_b^j;\boldsymbol{\theta})|^2.\]
Regarding the model architecture, for the sake of simplicity, we consider a neural network $u(x;\boldsymbol{\theta})$ with one hidden layer, activated by the function $\sigma=\max\{0,x^3\}$, which was considered in some previous literature (see, e.g., \cite{drm}):
\begin{equation}\label{NN_def}
u(x;\boldsymbol{\theta})=\frac{1}{\sqrt{d_1}}\sum_{k=1}^{d_1}W_2^k\sigma(W_1^kx+b_1^k)+b_2,
\end{equation}
where $d_1$ denotes the number of neurons in the hidden layer, and the parameters $\{W^k_1\}_{k=1}^{d_1}$, $\{W^k_2\}_{k=1}^{d_1}$, $\{b^k_1\}_{k=1}^{d_1}$ and $b_2$ are i.i.d. randomly chosen from the Gaussian distribution $\mathcal{N}(0,1)$. Before proceeding further, we introduce the following lemma which will be used to compute the averaged eigenvalues (see, e.g., \cite{math_sta}).
\begin{lemma}\label{Expectation}
    If $X\sim \mathcal{N}(0,\delta^2)$, then there holds
    \[\mathbb{E}\left[X^{2n}\right]=\frac{(2n)!}{2^nn!}\delta^{2n}\quad{\rm{and}}\quad \mathbb{E}\left[X^{2n-1}\right]=0 \qquad\forall n\in \mathbb{N}.\]
\end{lemma}

\subsection{Averaged Convergence Rate}\label{subsec:bdy_NTK}
In this section, we will explicitly compute the quantities ${\rm Tr}(\boldsymbol{K}_{uu}(0))/{N_b}$ and ${\rm Tr}(\boldsymbol{K}_{rr}(0))/{N_r}$ for the scaled problem \eqref{theory_scaled} on the given collocation points $\{x_b^j\}_{j=1}^{N_b}$ and $\{x_r^i\}_{i=1}^{N_r}$ with $x_b^j\sim\, {\rm{i.i.d.}}\,\,\mathcal{U}(\partial[0,N])$ and $x_r^i\sim\, {\rm{i.i.d.}}\,\,\mathcal{U}((0,N))$.
We first compute the convergence rate regarding the boundary collocation points. From the definition \eqref{NTK_PINN}, we see that
\[
\frac{{\rm Tr}(\boldsymbol{K}_{uu}(0))}{N_b}=\frac{1}{N_b}\sum_{j=1}^{N_b}\left\langle\frac{du(x_b^j;\boldsymbol{\theta})}{d\boldsymbol{\theta}},\frac{du(x_b^j;\boldsymbol{\theta})}{d\boldsymbol{\theta}}\right\rangle,
\]
where each summand can be decomposed as

\begin{align*}
\left\langle\frac{du(x_b^j;\boldsymbol{\theta})}{d\boldsymbol{\theta}},\frac{du(x_b^j;\boldsymbol{\theta})}{d\boldsymbol{\theta}}\right\rangle
&=\left<\frac{du(x_b^j;\boldsymbol{\theta})}{dW_1},\frac{du(x_b^j;\boldsymbol{\theta})}{dW_1}\right>+\left<\frac{du(x_b^j;\boldsymbol{\theta})}{dW_2},\frac{du(x_b^j;\boldsymbol{\theta})}{dW_2}\right>\\
&\hspace{4mm}+\left<\frac{du(x_b^j;\boldsymbol{\theta})}{db_1},\frac{du(x_b^j;\boldsymbol{\theta})}{db_1}\right>+\left<\frac{du(x_b^j;\boldsymbol{\theta})}{db_2},\frac{du(x_b^j;\boldsymbol{\theta})}{db_2}\right>.
\end{align*}
Since the parameters $W_1^k$, $W_2^k$ and $b_1^k$ are i.i.d. random samples from $\mathcal{N}(0,1)$ for each $k=1,\cdots d_1$, we may write them as random variables without the index $k$ 
\[
X_1:=W^k_1\sim\,\mathcal{N}(0,1),\,\, X_2:=W^k_2\sim\,\mathcal{N}(0,1),\,\,Z:=b^k_1\sim\,\mathcal{N}(0,1).
\]
Furthermore, since the random variables $X_1$ and $Z$ are independent, it is easy to verify that
\[
Y:=X_1x^j_b+Z\sim\mathcal{N}(0,|x_b^j|^2+1).
\]
Then if the number of neurons in the hidden layer  $d_1$ goes to infinity, we have from the law of large numbers and Lemma \ref{Expectation} that
\begin{equation*}
\begin{aligned}
\left<\frac{du(x_b^j;\boldsymbol{\theta})}{dW_1},\frac{du(x_b^j;\boldsymbol{\theta})}{dW_1}\right> 
&= \frac{1}{d_1}\sum_{k=1}^{d_1}\left|W_2^k 
\sigma'(W_1^k x_b^j+b_1^k)x_b^j\right|^2\xrightarrow{\mathcal{P}}|x_b^j|^2\mathbb{E}\left[\sigma'(Y)^2\right],\\
\left<\frac{du(x_b^j;\boldsymbol{\theta})}{dW_2},\frac{du(x_b^j;\boldsymbol{\theta})}{dW_2}\right> 
&= \frac{1}{d_1}\sum_{k=1}^{d_1}\left|\sigma(W_1^k x_b^j+b_1^k)\right|^2 \xrightarrow{\mathcal{P}} \mathbb{E}\left[\sigma(Y)^2\right],\\
\left<\frac{du(x_b^j;\boldsymbol{\theta})}{db_1},\frac{du(x_b^j;\boldsymbol{\theta})}{db_1}\right> 
&= \frac{1}{d_1}\sum_{k=1}^{d_1}\left|W_2^k 
\sigma'(W_1^k x_b^j+b_1^k)\right|^2\xrightarrow{\mathcal{P}}\mathbb{E}\left[\sigma'(Y)^2\right],\\
\left<\frac{du(x_b^j;\boldsymbol{\theta})}{db_2},\frac{du(x_b^j;\boldsymbol{\theta})}{db_2}\right> &= 1.\\
\end{aligned}
\end{equation*}
By recalling the definition of the activation function $\sigma=\max\{0,x^3\}$ and using Lemma \ref{Expectation}, there holds
\begin{align*}
    \mathbb{E}\left[\sigma'(Y)^2\right]
    &=\frac{9}{2\sqrt{2\pi(|x^j_b|^2+1)}}\int^{\infty}_{-\infty}y^4\exp\left(-\frac{y^2}{2(|x^j_b|^2+1)}\right)\dy=\frac{9}{2}\mathbb{E}\left[Y^4\right]=\frac{27}{2}\left(|x^j_b|^2+1\right)^2.
    \end{align*}
Next, by Lemma \ref{Expectation} again, the second term can be computed as
\begin{align*}
    \mathbb{E}\left[\sigma(Y)^2\right]
    &=\frac{1}{2\sqrt{2\pi(|x^j_b|^2+1)}}\int^{\infty}_{-\infty}y^6\exp\left(-\frac{y^2}{2(|x^j_b|^2+1)}\right)\dy=\frac{1}{2}\mathbb{E}\left[Y^6\right]=\frac{15}{2}\left(|x^j_b|^2+1\right)^3.
    \end{align*}
Altogether, we obtain as $d_1\rightarrow\infty$ that
\begin{equation}\label{bdry_tr}
\begin{aligned}
\frac{1}{N_b}\sum_{j=1}^{N_b}\left<\frac{du(x_b^j;\boldsymbol{\theta})}{d\boldsymbol{\theta}},\frac{du(x_b^j;\boldsymbol{\theta})}{d\boldsymbol{\theta}}\right>
&\xrightarrow{\mathcal{P}}\frac{1}{N_b}\sum_{j=1}^{N_b}\left(21|x^j_b|^6+63|x^j_b|^4+63|x^j_b|^2+21\right).
\end{aligned}
\end{equation}
Therefore, for sufficiently large $N_b$, by the law of large numbers and the fact $x^j_b\sim\mathcal{U}(\partial[0,N])$, we can conclude that

\begin{equation}\label{data_tr_final}
\frac{{\rm Tr}(\boldsymbol{K}_{uu}(0))}{N_b}=\frac{1}{N_b}\sum_{j=1}^{N_b}\left\langle\frac{du(x_b^j;\boldsymbol{\theta})}{d\boldsymbol{\theta}},\frac{du(x_b^j;\boldsymbol{\theta})}{d\boldsymbol{\theta}}\right\rangle\approx \mathcal{O}(N^6).
\end{equation}

As a next step, we shall compute the value ${\rm Tr}(\boldsymbol{K}_{\rm rr}(0))$ on the data set $\{x_r^i\}_{i=1}^{N_r}$ with $x^i_r\sim\,{\rm{i.i.d.}}\,\,\mathcal{U}((0,N))$. By recalling the definition \eqref{NTK_PINN} and \eqref{NN_def}, we see that

\begin{equation}\label{res_trace}
    \frac{{\rm Tr}(\boldsymbol{K}_{rr}(0))}{N_r}=\frac{1}{N_r}\sum_{i=1}^{N_r}\left<\frac{du_{xx}(x_r^i;\boldsymbol{\theta})}{d\boldsymbol{\theta}},\frac{du_{xx}(x_r^i;\boldsymbol{\theta})}{d\boldsymbol{\theta}}\right>
\end{equation}
where we can verify that
\begin{equation}\label{xx_wave}
    u_{xx}(x_r^i;\boldsymbol{\theta})=\frac{1}{\sqrt{d_1}}\sum_{k=1}^{d_1}W_2^k\sigma''(W_1^k x_r^i + b_1^k)|W_1^k|^2.
\end{equation}
In \eqref{res_trace}, for each $i=1,\cdots,N_r$, the summand can be written as
\begin{align*}
    \left<\frac{du_{xx}(x_r^i;\boldsymbol{\theta})}{d\boldsymbol{\theta}},\frac{du_{xx}(x_r^i;\boldsymbol{\theta})}{d\boldsymbol{\theta}}\right>
    &=\left<\frac{du_{xx}(x_r^i;\boldsymbol{\theta})}{dW_1},\frac{du_{xx}(x_r^i;\boldsymbol{\theta})}{dW_1}\right>+\left<\frac{du_{xx}(x_r^i;\boldsymbol{\theta})}{dW_2},\frac{du_{xx}(x_r^i;\boldsymbol{\theta})}{dW_2}\right>\\
    &\hspace{4mm}+\left<\frac{du_{xx}(x_r^i;\boldsymbol{\theta})}{db_1},\frac{du_{xx}(x_r^i;\boldsymbol{\theta})}{db_1}\right>+\left<\frac{du_{xx}(x_r^i;\boldsymbol{\theta})}{db_2},\frac{du_{xx}(x_r^i;\boldsymbol{\theta})}{db_2}\right>.
\end{align*}
For this case, we should be more careful since $X_1$ and $X_1x^i_r+Z$ appear simultaneously in \eqref{xx_wave} and they may not be independent. As $d_1 \rightarrow \infty$, by Lemma \ref{Expectation} and the law of large numbers, we obtain
\begin{align*}
\left<\frac{du_{xx}(x_r^i;\boldsymbol{\theta})}{dW_1},\frac{du_{xx}(x_r^i;\boldsymbol{\theta})}{dW_1}\right> &= \frac{1}{d_1}\sum_{k=1}^{d_1}\left||W_1^k|^2W_2^k 
\sigma^{(3)}(W_1^k x_r^i+b_1^k)x_r^i+2W_1^kW_2^k\sigma''(W_1^kx_r^i+b_1^k)\right|^2\\
&\xrightarrow{\mathcal{P}}|x^i_r|^2\mathbb{E}\left[X_1^4\sigma^{(3)}(X_1x^i_r+Z)^2\right]+4x_r^i\mathbb{E}\left[X_1^3\sigma^{(3)}(X_1x^i_r+Z)\sigma''(X_1x^i_r+Z)\right]\\
&\hspace{6mm}+4\mathbb{E}\left[X_1^2\sigma''(X_1x^i_r+Z)^2\right],\\
\left<\frac{du_{xx}(x_r^i;\boldsymbol{\theta})}{dW_2},\frac{du_{xx}(x_r^i;\boldsymbol{\theta})}{dW_2}\right> &= \frac{1}{d_1}\sum_{k=1}^{d_1}\left|\sigma''(W_1^k x_r^i+b_1^k)|W_1^k|^2\right|^2\xrightarrow{\mathcal{P}}\mathbb{E}\left[X_1^4\sigma''(X_1x^i_r+Z)^2\right],\\
\left<\frac{du_{xx}(x_r^i;\boldsymbol{\theta})}{db_1},\frac{du_{xx}(x_r^i;\boldsymbol{\theta})}{db_1}\right> &= \frac{1}{d_1}\sum_{k=1}^{d_1}\left|W_2^k 
\sigma^{(3)}(W_1^k x_r^i+b_1^k)|W_1^k|^2\right|^2\xrightarrow{\mathcal{P}}\mathbb{E}\left[X_1^4\sigma^{(3)}(X_1x^i_r+Z)^2\right],\\
\left<\frac{du_{xx}(x_r^i;\boldsymbol{\theta})}{db_2},\frac{du_{xx}(x_r^i;\boldsymbol{\theta})}{db_2}\right> &=0.\\
\end{align*}
Regarding the derivative of $\sigma$ in the sense of distribution to define $\sigma^{(3)}$, we see by Lemma \ref{Expectation} and the fact $|\sigma^{(3)}(X_1x^i_r+Z)|\leq 6$ that
\[\mathbb{E}\left[X_1^4\sigma^{(3)}(X_1x^i_r+Z)^2\right]\leq\frac{36}{{2\pi}} \int_{-\infty}^{\infty}\int_{-\infty}^{\infty} x^4 \exp\left(-\frac{x^2+z^2}{2}\right)\dx\dz=36\mathbb{E}\left[X_1^4\right]=108.\]
Next, by Lemma \ref{Expectation} and Fubini's theorem,
\begin{align*}
\mathbb{E}\left[X_1^3\sigma^{(3)}(X_1x^i_r+Z)\sigma''(X_1x^i_r+Z)\right]
&\leq\frac{6}{2\pi}\int^{\infty}_{-\infty}\int^{\infty}_{-\infty}x^3\sigma''(xx^i_r+z)\exp\left(-\frac{x^2+z^2}{2}\right)\dz\dx\\
&\leq\frac{36}{2\pi} \int_{-\infty}^{\infty}\int_{-x^i_rx}^{\infty} (x_r^ix^4+x^3z) \exp\left(-\frac{x^2+z^2}{2}\right)\dz\dx\\
&\leq\frac{36}{2\pi}\int^{\infty}_{-\infty}\int^{\infty}_{-\infty}x^i_rx^4\exp\left(-\frac{x^2+z^2}{2}\right)\dz\dx\\
&\hspace{5mm}+\frac{36}{2\pi}\int^{\infty}_{-\infty}x^3\exp\left(-\frac{x^2}{2}\right)\int_{-x^i_rx}^{\infty}z\exp\left(-\frac{z^2}{2}\right)\dz\dx\\
&= 36x_r^i\mathbb{E}\left[X_1^4\right]+\frac{36}{2\pi}\int^{\infty}_{-\infty}x^3\exp\left(-\frac{(|x^i_r|^2+1)x^2}{2}\right)=108x^i_r.
\end{align*}
Moreover, by Lemma \ref{Expectation} again and the fact $|\sigma''(X_1x^i_r+Z)|^2\leq 36|X_1x^i_r+Z|^2$, we have
\begin{align*}
    \mathbb{E}\left[X_1^{2}\sigma''(X_1x^i_r+Z)^2\right]
    &\leq\frac{36}{2\pi}\int^{\infty}_{-\infty}\int^{\infty}_{-\infty}x^{2}(x^2|x_r^i|^2+2xzx_r^i+z^2)\exp\left(-\frac{x^2+z^2}{2}\right)\dz\dx\\   &=36|x_r^i|^2\mathbb{E}\left[X_1^{4}\right]+36\mathbb{E}\left[X_1^{2}\right]\mathbb{E}\left[Z^2\right]=36\left(3|x^i_r|^2+1\right),
\end{align*}
and similarly we obtain
\begin{align*}
    \mathbb{E}\left[X_1^{4}\sigma''(X_1x^i_r+Z)^2\right]\leq36|x_r^i|^2\mathbb{E}\left[X_1^{6}\right]+36\mathbb{E}\left[X_1^{4}\right]\mathbb{E}\left[Z^2\right]=36\left(15|x^i_r|^2+3\right).
\end{align*}
From above, by collecting all the terms, we have as $d_1\rightarrow\infty$ that
\[
    \frac{1}{N_r}\sum_{i=1}^{N_r}\left<\frac{du_{xx}(x_r^i;\boldsymbol{\theta})}{d\boldsymbol{\theta}},\frac{du_{xx}(x_r^i;\boldsymbol{\theta})}{d\boldsymbol{\theta}}\right>
    \lesssim\frac{1}{N_r}\sum_{i=1}^{N_r}\left(|x^i_r|^2+1\right).
\]
Finally, for sufficiently large $N_r$, from the law of large numbers and the fact $x^i_r\sim\mathcal{U}((0,N))$, we may conclude that
\begin{equation}\label{res_tr_final}
\frac{{\rm Tr}(\boldsymbol{K}_{rr}(0))}{N_r}=\frac{1}{N_r}\sum_{j=1}^{N_r}\left\langle\frac{du_{xx}(x_r^i;\boldsymbol{\theta})}{d\boldsymbol{\theta}},\frac{du_{xx}(x_r^i;\boldsymbol{\theta})}{d\boldsymbol{\theta}}\right\rangle\lesssim \mathcal{O}(N^2).
\end{equation}

\begin{figure}
    \centering
    \includegraphics[width=0.9\textwidth]{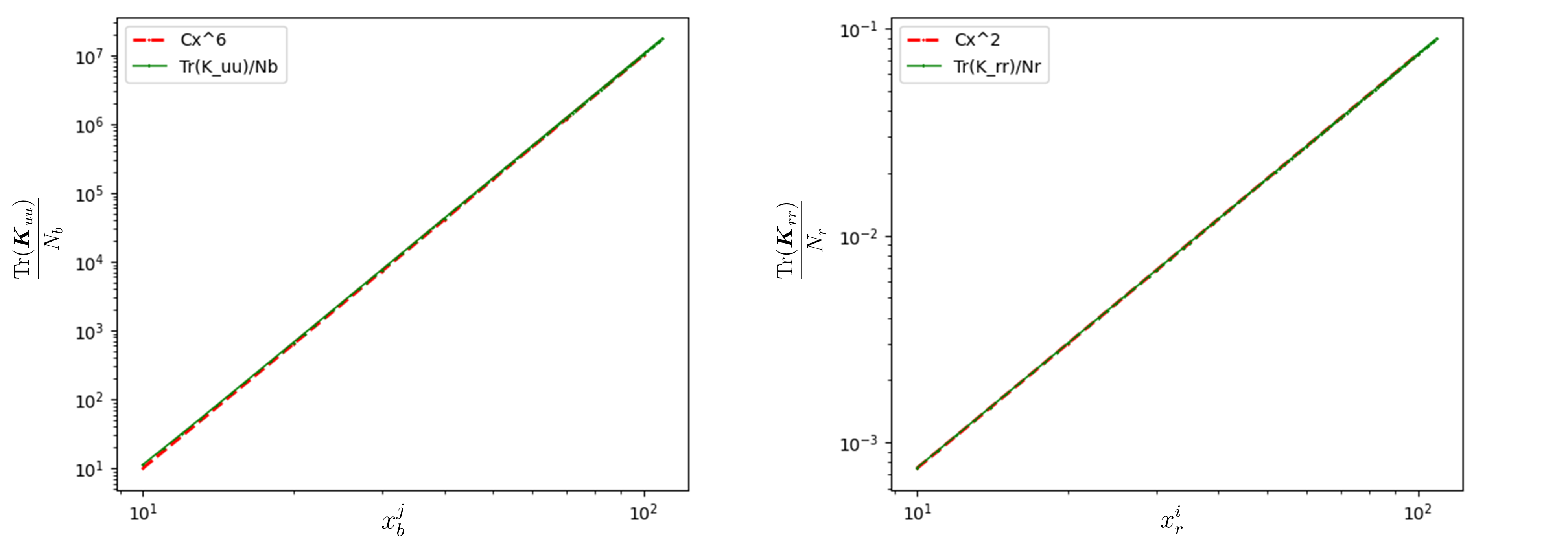}
    \caption{Log-log plots of $\frac{{\rm Tr}(\boldsymbol{K}_{uu})}{N_b}$ and $\frac{{\rm Tr}(\boldsymbol{K}_{rr})}{N_r}$ against the input samples $x^j_b$ and $x^i_r$. We can confirm that the experimental results align well with the theoretical results.}
    \label{log_log}
\end{figure}
\begin{figure}
    \centering
    \includegraphics[width=0.9\textwidth]{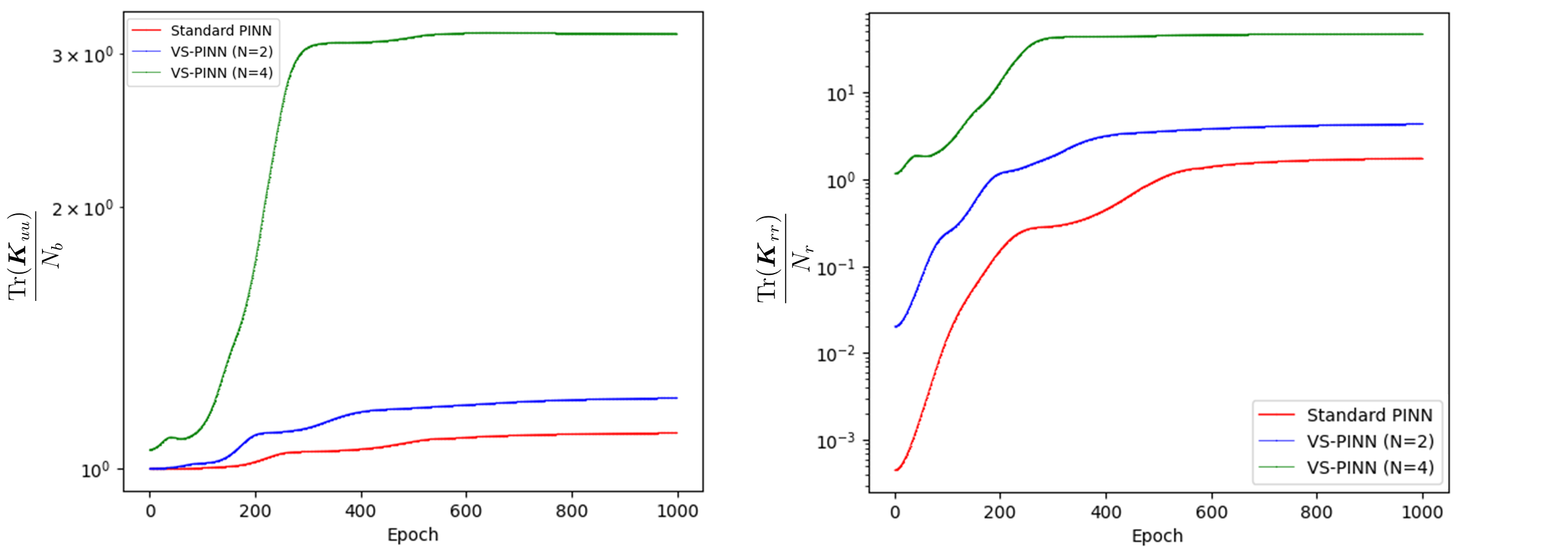}
    \caption{The averaged convergence rates $\frac{{\rm Tr}(\boldsymbol{K}_{uu})}{N_b}$ and $\frac{{\rm Tr}(\boldsymbol{K}_{rr})}{N_r}$ during training.}
    \label{cubic_graph}
\end{figure}

Therefore, from \eqref{data_tr_final} and \eqref{res_tr_final}, we have shown that the averaged convergence rate \eqref{avg_conv_rate} becomes larger when we apply the variable scaling. Figure \ref{log_log} displays the averaged convergence rates $\frac{{\rm Tr}(\boldsymbol{K}_{uu})}{N_b}$ and $\frac{{\rm Tr}(\boldsymbol{K}_{rr})}{N_r}$
as functions of the input samples $x^j_b$ and $x^i_r$ respectively. From the figures, we can confirm that the experimental results go along with the theoretical findings. Furthermore, from Figure \ref{cubic_graph}, we can observe that the convergence rates of the VS-PINN are maintained higher as N increases during the training process, which supports the faster convergence of the VS-PINN observed in the experiments. Note, however, that this result only deals with the theoretical convergence rates and holds for the ideal case. In practice, large scaling factor $N$ does not always guarantee faster convergence, and the choice of suitable $N$ is important. We will discuss this issue in the following section.

\subsection{Proper choice of the scaling factor $N$}\label{subsec:NTK_Poisson}
In the previous sections, we have seen that applying the variable scaling to PINNs significantly improves the performance, especially for PDE problems with stiff behavior and high frequency. Furthermore, throughout the NTK analysis, we have provided some theoretical evidence of why the VS-PINN can achieve better performance as the scaling factor $N$ increases. However, simply increasing $N$ without careful consideration may deteriorate the performance of the VS-PINN in practice. We shall briefly discuss this issue in this section through a simple experiment.

Let us consider the one-dimensional Poisson
problem
\begin{equation*}
    \begin{aligned}
        -u_{xx} &= \pi^2\sin(\pi x)\;\;\;\; x\in[0,1],\\
        u(0)&=u(1)=0,
    \end{aligned}
\end{equation*}
and compute the approximate solution using the VS-PINN with various scaling factors. We shall use a two-layer neural network with $40000$ neurons in the hidden layer, where the weight and bias parameters are initialized according to the normal distribution $\mathcal{N}(0,0.1^2)$. The training proceeds with two boundary data points and $50$ uniformly-distributed interior data points.

\begin{figure}
    \centering
    \includegraphics[width=0.9\textwidth]{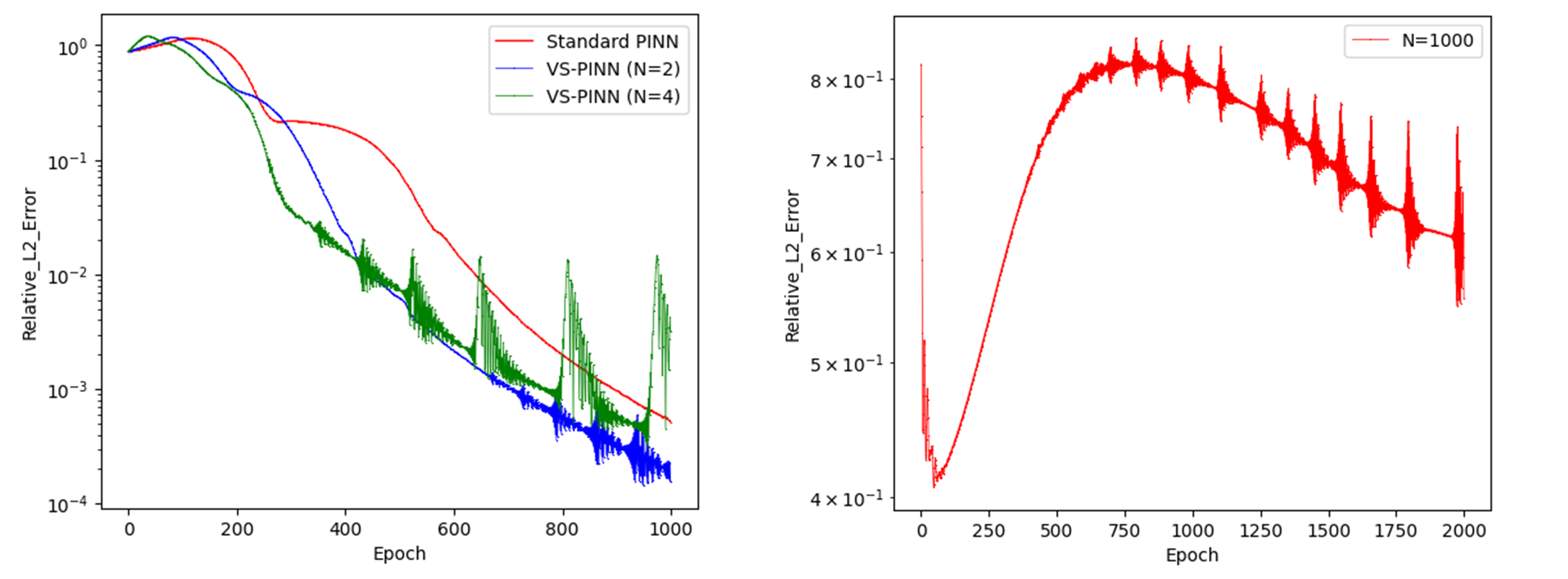}
    \caption{Learning curves of the VS-PINN with various scaling factors.}
    \label{VS_fail}
\end{figure}

As we can see from Figure \ref{VS_fail}, when $N=2$, the performance of the VS-PINN is better than the standard PINN. However, if we set $N=4$, we can observe that the training becomes unstable and accuracy decreases after a certain point. Furthermore, if we choose a much larger scaling factor $N=1000$, the performance deteriorates even further and the VS-PINN is no longer promising compared to the original PINN. In fact, when the scaling factor is excessively large, the training dynamic becomes very sensitive and we need to carefully adjust the parameters related to the training. For example, since the domain is getting extended as $N$ increases, we need a larger number of collocation points to guarantee satisfactory performance. This observation aligns with the situations encountered when applying the law of large numbers in the theoretical analysis in the previous section. In other words, if the scaling factor is chosen appropriately, we can effectively train the model only with a reasonable number of sample points leveraging the advantages of the VS-PINN, while this is no longer the case for the extreme choice of $N$. Moreover, if $N$ is too large, the corresponding scaled loss function might become excessively large. In this case, the change of parameters in a single step of gradient descent becomes larger, which hinders us from finding the global minima and makes the training unstable. Therefore, as $N$ increases, we need to put in a lot of effort to find the proper settings for training and need to carefully adjust the factors related to training such as the weight parameters in the loss, variance of the initialization distribution and the learning rate. In conclusion, to effectively utilize VS-PINN without excessive effort to find a proper circumstance for training, it is important to choose a reasonable scaling factor $N$, which can indeed significantly improve the performance of the PINNs as we saw in the experiments in Section \ref{sec:exp}.






\section{Conclusion}\label{sec:conclude}
In this paper, we proposed a novel training method for PINNs based on the variable-scaling technique. This method is simple and it can be applied to a wide range of PDE problems with stiff behavior or high frequency. We applied this new strategy to various problems including the wave equation, the Allen-Cahn equation, the singular perturbation problem and the Navier--Stokes equations. Throughout these numerical experiments, we confirmed that the VS-PINN provides a significant performance improvement both in accuracy and training efficiency for stiff problems and even for nonlinear problem. Moreover, we performed the NTK analysis to provide theoretical evidence of why the VS-PINN works well. More precisely, we explicitly computed the eigenvalues of the NTK of VS-PINNs which concerns the rate of convergence for training, and showed that the eigenvalues corresponding to the VS-PINN are indeed larger when the variable is scaled.

An interesting future research direction is to theoretically analyze the convergence properties of the proposed method from the viewpoint of approximation and generalization. Furthermore, the systematic guidelines for choosing the scaling factor $N$ and the parameters related to training should also be investigated. Since increasing $N$ causes instability in training, establishing a general strategy that can address this issue while further enhancing the overall performance of the VS-PINN is an important research question. Since the choice of the parameters depends on the structure of the equation under consideration, we will provide a more detailed answer to this question with a general principle and theoretical intuition by examining a wider variety of PDE problems using the VS-PINN, which will be addressed in the forthcoming paper.

\vspace{5mm}

\bibliographystyle{abbrv}
\bibliography{sample}

\end{document}